\documentclass[a4paper,10pt]{article}
\usepackage{stmaryrd}
\usepackage{amsfonts}
\usepackage{bbm}
\usepackage{amscd}
\usepackage{mathrsfs}
\usepackage{latexsym,amssymb,amsmath,amscd,amscd,amsthm,amsxtra}
\usepackage[dvips]{graphicx}
\usepackage[utf8]{inputenc}
\usepackage[T1]{fontenc}
\usepackage{lmodern}
\usepackage{amssymb}
\usepackage[all]{xy}
\usepackage{nicefrac,mathtools,enumitem}
\usepackage{microtype}

\newtheorem{thm}{Theorem}[section]
\newtheorem{lem}[thm]{Lemma}

\newtheorem{pro}[thm]{Proposition}
\newtheorem{ex}[thm]{Example}
\newtheorem{rmk}[thm]{Remark}
\newtheorem{defi}[thm]{Definition}

\newcommand {\emptycomment}[1]{} 

\usepackage{color}

\newcommand{\be }{\begin{equation}}
\newcommand{\ee }{\end{equation}}

\newcommand{\Diff}{\mathsf {Diff}}
\newcommand{\DiffCat}{\mathsf {DiffCat}}

\newcommand{\pf}{\noindent{\bf Proof.}\ }

\newcommand{\curv}{\mathsf{curv}}
\newcommand{\fcurv}{\mathsf{fcurv}}
\newcommand{\h}{\mathbbm h}
\newcommand{\hh}{\mathrm h}

\newcommand{\HH}{\mathbb H}

\newcommand{\huaS}{\mathcal{S}}
\newcommand{\huaA}{\mathcal{A}}

\newcommand{\huaE}{\mathcal{E}}

\newcommand{\huaG}{\mathcal{G}}

\newcommand{\huaV}{\mathcal{V}}

\newcommand{\huaP}{\mathcal{P}}
\newcommand{\huaD}{\mathcal{D}}

\newcommand{\CWM}{C^{\infty}(M)}

\newcommand{\frkd}{\mathfrak d}

\newcommand{\frkf}{\mathfrak f}

\newcommand{\frkh}{\mathfrak h}

\newcommand{\frkk}{\mathfrak k}
\newcommand{\frkl}{\mathfrak l}

\newcommand{\frkt}{\mathfrak t}

\newcommand{\frkG}{\mathfrak G}

\newcommand{\frkX}{\mathfrak X}

\def\qed{\hfill ~\vrule height6pt width6pt depth0pt}

\newcommand{\half}{\frac{1}{2}}

\newcommand{\br}[1]{   [ \cdot,    \cdot  ]   }

\newcommand{\dev}{\mathfrak{D}}

\newcommand{\id}{\rm{id}}

\newcommand{\g}{\mathfrak g}

\newcommand{\ii}{\mathbbm{i}}
\newcommand{\jd}{\mathbbm{j}}

\newcommand{\dM}{\mathrm{d}}

\newcommand{\LWX}{\mathrm{CLWX}}
\newcommand{\Hom}{\mathrm{Hom}}
\newcommand{\Der}{\mathrm{Der}}
\newcommand{\DER}{\mathrm{DER}}

\newcommand{\Ad}{\mathrm{Ad}}

\newcommand{\gl}{\mathfrak {gl}}

\newcommand{\ad}{\mathrm{ad}}

\newcommand{\pr}{\mathrm{pr}}

\newcommand{\ve}{\mathrm{v}}

\newcommand{\sgn}{\mathrm{sgn}}

\begin{document}
\title{ The first Pontryagin class of a quadratic Lie 2-algebroid
 }
\author{Yunhe Sheng\\
Department of Mathematics, Jilin University,
Changchun 130012, China
\\ Email: shengyh@jlu.edu.cn
}

\date{}
\maketitle
\begin{abstract}
In this paper, first we give a detailed study on the structure of a transitive Lie 2-algebroid and describe a transitive Lie 2-algebroid using a morphism from the tangent Lie algebroid $TM$ to a strict  Lie 3-algebroid constructed from derivations. Then we introduce the notion of a quadratic Lie 2-algebroid and define its first Pontryagin class, which is a cohomology class in $H^5(M)$. Associated to a $\LWX$ $2$-algebroid, there is a quadratic Lie 2-algebroid naturally. Conversely, we show that the first Pontryagin class of a quadratic Lie 2-algebroid is the obstruction class of the existence of a $\LWX$-extension. Finally we construct a quadratic Lie 2-algebroid from a trivial principle  $2$-bundle with a $\Gamma$-connection and show that its first Pontryagin class is trivial.
\end{abstract}
\tableofcontents

\section{Introduction}

The notion of a Lie algebroid was introduced by Pradines in 1967, which
is a generalization of Lie algebras and tangent bundles. Just as Lie algebras are the
infinitesimal objects of Lie groups, Lie algebroids are the infinitesimal objects of Lie groupoids. See \cite{General theory of Lie groupoid and Lie algebroid} for general theory about Lie algebroids. A transitive Lie algebroid (i.e. $\rho$ is surjective) $(A,[\cdot,\cdot]_A,\rho)$ is called a quadratic Lie algebroid if there is a nondegenerate symmetric invariant bilinear form $\huaS$ on $\ker(\rho)$. A connection of a quadratic Lie algebroid is a splitting $\sigma:TM\longrightarrow A$ of the short exact sequence
$$
0\longrightarrow \ker(\rho)\longrightarrow A \stackrel{\rho}{\longrightarrow} TM\longrightarrow 0.
$$
The curvature $R\in\Hom(\wedge^2TM,\ker(\rho))$ is defined by
$$
R(X,Y)=\sigma[X,Y]-[\sigma(X),\sigma(Y)]_A,\quad \forall X,Y\in\Gamma(TM).
$$
Then the first Pontryagin class of a quadratic Lie algebroid is defined to be the cohomology class $[\huaS(R,R)]$ in $H^4(M)$.

The notion of a Courant algebroid was introduced in \cite{lwx} in the study of the double of a Lie bialgebroid. An alternative definition was given in \cite{Roytenbergphdthesis}. See the review article \cite{Schwarzbach4} for more information. Roughly speaking, a Courant algebroid is a vector bundle $E$, whose section space is a Leibniz algebra, together with an anchor map $\rho:E\longrightarrow TM$ and a nondegenerate symmetric bilinear form $S$, such that some compatibility conditions are satisfied. For a transitive Courant algebroid (i.e. $\rho$ is surjective), $E/\rho^*(T^*M)$ is a quadratic Lie algebroid. Conversely, it is proved in \cite{Bressler:Pclass,ChenRCA} that a quadratic Lie algebroid admits a Courant-extension if and only if its first Pontryagin class is trivial.

Recently, people have paid more attention to higher categorical
structures by reasons in both mathematics and physics. A Lie $2$-algebra
is the categorification of a Lie algebra \cite{baez:2algebras}. The Jacobi identity is
replaced by a natural transformation, called the Jacobiator, which
also satisfies some coherence laws of its own. If a skew-symmetric bracket is used in the definition of a Courant algebroid,  the authors showed that the underlying algebraic structure of a Courant algebroid is a Lie 2-algebra  \cite{rw}. Along this approach, the author defined split Lie $n$-algebroids in \cite{sz} using the language of graded vector bundles. Usually an NQ-manifold of degree $n$ is considered as a Lie $n$-algebroid \cite{Voronov:2010halgd}. The equivalence between  the category of  split Lie $n$-algebroids and the category of NQ-manifolds of degree $n$ was given in \cite{BP}.

$\LWX$ 2-algebroids (named after Courant-Liu-Weinstein-Xu) were introduced in \cite{LiuSheng}, which can be viewed  as the categorification of Courant algebroids, and one-to-one correspond to symplectic NQ-manifolds of degree 3. See \cite{Ikeda} for more applications in 4D topological field theory. Similar as the case of Courant algebroids, we can obtain a quadratic Lie 2-algebroid (called the ample Lie 2-algebroid) from an exact $\LWX$ 2-algebroid by modulo $\rho^*(T^*M)$. Then it is natural to ask the following question:
\begin{itemize}
   \item[$\bullet$]Whether every quadratic Lie 2-algebroid admits a $\LWX$-extension? If not, what is the obstruction?
\end{itemize}

 Motivated by the above question, first we study the structure of a transitive Lie 2-algebroid $(A_{-1}\oplus A_0;\rho,l_1,l_2,l_3)$. It turns out that a transitive Lie 2-algebroid is a nonabelian extension of the tangent Lie algebroid $TM$ by a graded bundle of Lie 2-algebra $A_{-1}\oplus\ker(\rho)$. By choosing a splitting (a connection), we give precise formulas of a transitive Lie 2-algebroid. One important property that we need to stress here is that the corresponding curvatures include a $\ker(\rho)$-valued 2-form and an $A_{-1}$-valued 3-form. Furthermore, we construct a strict Lie 3-algebroid using derivations of $A_{-1}\oplus\ker(\rho)$ and characterize a transitive Lie 2-algebroid by a morphism from the tangent Lie algebroid $TM$ to this strict Lie 3-algebroid. Then we go on to give the notion of a quadratic Lie 2-algebroid. For a quadratic Lie 2-algebroid,  we construct a 5-form using the aforementioned curvatures. We show that this 5-form is closed and does not depend on the choices of connections. Thus, it gives rise to a characteristic class, which we call the  first Pontryagin class of a quadratic Lie 2-algebroid. The first Pontryagin class of a quadratic Lie 2-algebroid plays the same role as the one for a quadratic Lie algebroid, namely, a quadratic Lie 2-algebroid admits a $\LWX$-extension if and only if its first Pontryagin class is trivial. This is the reason why we call this characteristic class the first Pontryagin class of a quadratic Lie 2-algebroid.

 In the past ten years, principle 2-bundles and higher gauge theory are deeply studied in \cite{ACJ,baez-schreiber,breen-messing,GP,MM,NSS,schreiber,schommer,waldorf:global,waldorf:parallel,woc11}.
 For a trivial principle $\Gamma$-$2$-bundle with a connection, where $\Gamma$ is a strict Lie 2-group, we construct a transitive Lie 2-algebroid which can be viewed as its Atiyah algebroid. We further show that if the Lie 2-algebra corresponding to the Lie 2-group $\Gamma$ is quadratic, then this Lie 2-algebroid is quadratic and its first Pontryagin class is trivial. We will go on to study the infinitesimal of principle 2-bundles in a separate paper.

 The paper is organized as follows. In Section 2, we recall Leibniz 2-algebras, Lie $n$-algebroids and $\LWX$ 2-algebroids. In Section 3,   we study transitive Lie 2-algebroids. In Subsection 3.1, we give precise conditions that the structure maps of a transitive Lie 2-algebroid satisfy.  In Subsection 3.2, we construct a strict Lie 3-algebroid using derivations and  characterize a transitive Lie 2-algebroid by a morphism from the tangent Lie algebroid $TM$ to this strict Lie 3-algebroid. In Section 4, we give the notion of a quadratic Lie 2-algebroid and define its first Pontryagin class, which is a cohomology class in $H^5(M)$. In Section 5, we study exact $\LWX$ 2-algebroids. In Subsection 5.1, we show that skeletal exact $\LWX$ 2-algebroids are classified by the higher analogue of the \v{S}evera class, which is a cohomology class in $H^4(M)$. In Subsection 5.2, we show that a $\LWX$ 2-algebroid gives rise to a quadratic Lie 2-algebroid and conversely, a quadratic Lie 2-algebroid admits a $\LWX$-extension if and only if its first Pontryagin class is trivial. In Section 6, for a trivial principle 2-bundle with a $\Gamma$-connection, we construct a transitive Lie 2-algebroid and show that its first Pontryagin class is trivial.

 \vspace{2mm}
 \noindent {\bf Acknowledgement:} We give our warmest thanks to Zhangju Liu, Konrad Waldorf,  Xiaomeng Xu and Chenchang Zhu for very useful comments and discussions. We also give our special thanks to the referee for very helpful suggestions that improve the paper.
This research is supported by NSFC (11471139) and NSF of Jilin Province (20170101050JC).

\section{Preliminaries}
\subsection{Leibniz 2-algebras}

As a  model for ``Leibniz
algebras that satisfy Jacobi identity up to all higher homotopies'',
  the notion of a strongly homotopy
Leibniz algebra, or a $Lod_\infty$-algebra was given in \cite{livernet} by Livernet,
 which was further studied by Ammar and Poncin
in \cite{ammardefiLeibnizalgebra}.
  In \cite{Leibniz2al}, the authors introduced the notion of a Leibniz 2-algebra, which is the categorification of a Leibniz algebra, and proved that the category of Leibniz 2-algebras and the category of 2-term $Lod_\infty$-algebras are equivalent.

\begin{defi}\label{defi:2leibniz}
  A   {\bf Leibniz $2$-algebra} $\huaV$ consists of the following data:
\begin{itemize}
\item[$\bullet$] a complex of vector spaces $\huaV:V_{-1}\stackrel{\dM}{\longrightarrow}V_0,$

\item[$\bullet$] bilinear maps $l_2:V_{-i}\times V_{-j}\longrightarrow
V_{-i-j}$, where  $0\leq i+j\leq1$,

\item[$\bullet$] a  trilinear map $l_3:V_0\times V_0\times V_0\longrightarrow
V_{-1}$,
   \end{itemize}
   such that for all $w,x,y,z\in V_0$ and $m,n\in V_{-1}$, the following equalities are satisfied:
\begin{itemize}
\item[$\rm(a)$] $\dM l_2(x,m)=l_2(x,\dM m),$
\item[$\rm(b)$]$\dM l_2(m,x)=l_2(\dM m,x),$
\item[$\rm(c)$]$l_2(\dM m,n)=l_2(m,\dM n),$
\item[$\rm(d)$]$\dM l_3(x,y,z)=l_2(x,l_2(y,z))-l_2(l_2(x,y),z)-l_2(y,l_2(x,z)),$
\item[$\rm(e_1)$]$ l_3(x,y,\dM m)=l_2(x,l_2(y,m))-l_2(l_2(x,y),m)-l_2(y,l_2(x,m)),$
\item[$\rm(e_2)$]$ l_3(x,\dM m,y)=l_2(x,l_2(m,y))-l_2(l_2(x,m),y)-l_2(m,l_2(x,y)),$
\item[$\rm(e_3)$]$ l_3(\dM m,x,y)=l_2(m,l_2(x,y))-l_2(l_2(m,x),y)-l_2(x,l_2(m,y)),$
\item[$\rm(f)$] the Jacobiator identity:\begin{eqnarray*}
&&l_2(w,l_3(x,y,z))-l_2(x,l_3(w,y,z))+l_2(y,l_3(w,x,z))+l_2(l_3(w,x,y),z)\\
&&-l_3(l_2(w,x),y,z)-l_3(x,l_2(w,y),z)-l_3(x,y,l_2(w,z))\\
&&+l_3(w,l_2(x,y),z)+l_3(w,y,l_2(x,z))-l_3(w,x,l_2(y,z))=0.\end{eqnarray*}
   \end{itemize}
\end{defi}
We usually denote a   Leibniz 2-algebra by
$(V_{-1}\oplus V_0;\dM,l_2,l_3)$, or simply by
$\huaV$. In particular, if both $l_2$ and $l_3$ are skew-symmetric, we obtain the notion of a Lie 2-algebra \cite{baez:2algebras}.

\begin{defi}\label{defi:Leibniz morphism}
 Let $\huaV$ and $\huaV^\prime$ be   Leibniz $2$-algebras, a morphism $\frkf$
 from $\huaV$ to $\huaV^\prime$ consists of
\begin{itemize}
  \item[$\bullet$] linear maps $f_0:V_0\longrightarrow V_0^\prime$
  and $f_1:V_1\longrightarrow V_1^\prime$ commuting with the
  differential, i.e.
  $$
f_0\circ \dM=\dM^\prime\circ f_1;
  $$
  \item[$\bullet$] a bilinear map $f_2:V_0\times V_0\longrightarrow
  V_1^\prime$,
\end{itemize}
  such that  for all $x,y,z\in L_0,~m\in
L_1$, we have
\begin{equation}\label{eqn:DGLA morphism c 1}\left\{\begin{array}{rll}
f_0l_2(x,y)-l_2^\prime(f_0(x),f_0(y))&=&\dM^\prime
f_2(x,y),\\
f_1l_2(x,m)-l_2^\prime(f_0(x),f_1(m))&=&f_2(x,\dM
m),\\
f_1l_2(m,x)-l_2^\prime(f_1(m),f_0(x))&=&f_2(\dM
m,x),\end{array}\right.
\end{equation}
and
\begin{eqnarray}
\nonumber&&-f_1(l_3(x,y,z))+l_2^\prime(f_0(x),f_2(y,z))-l_2^\prime(f_0(y),f_2(x,z))-l_2^\prime(f_2(x,y),f_0(z))\\
\label{eqn:DGLA morphism c 2}
&&-f_2(l_2(x,y),z)+f_2(x,l_2(y,z))-f_2(y,l_2(x,z))+l_3^\prime(f_0(x),f_0(y),f_0(z))=0.
\end{eqnarray}
\end{defi}

\subsection{Lie $n$-algebroids and $\LWX$ 2-algebroids}

The notion of a split Lie $n$-algebroid was introduced in \cite{sz}.

\begin{defi}
A {\bf split Lie $n$-algebroid} is a graded vector bundle $\huaA= A_{-n+1}\oplus \cdots \oplus A_{-1}\oplus A_0$ over a manifold $M$ equipped with a bundle map (the anchor) $\rho:A_0\longrightarrow TM$, and  brackets $l_i:\Gamma(\wedge^i\huaA)\longrightarrow \Gamma(\huaA)$ with degree $2-i$ for $i=1,2,\cdots, l_{n+1}$, such that
\begin{itemize}
\item[$\rm(1)$]$(\Gamma(\huaA);l_1,l_2,\cdots, l_{n+1})$ is a Lie $n$-algebra ($n$-term $L_\infty$-algebra);
\item[$\rm(2)$]$l_2$ satisfies the Leibniz rule with respect to the anchor $\rho$:
$$l_2(X^0,fY)=fl_2(X^0,Y)+\rho(X^0)(f)Y,\quad \forall X^0\in\Gamma(A_0),f\in\CWM,Y\in\Gamma(\huaA);$$
\item[$\rm(3)$]for $i\neq2$, $l_i$ are $\CWM$-linear.
 \end{itemize}
 Denote a Lie $n$-algebroid  by $(\huaA;\rho,l_1,l_2,\cdots, l_{n+1}).$
\end{defi}
A split Lie $n$-algebroid is said to be {\bf strict} if $l_i=0$ for all $i>2$.
In this paper, we will only use split Lie $2$-algebroids and strict split Lie $3$-algebroids. See \cite{BP}
for more details about the category of Lie $n$-algebroids.

\begin{lem}\label{lem:propertyLien}
 Let $(\huaA;\rho,l_1,l_2,\cdots, l_{n+1})$ be a Lie $n$-algebroid. Then we have
 \begin{eqnarray}
   \rho\circ l_1&=&0,\\
  \rho l_2(X^0,Y^0)&=&[\rho(X^0),\rho(Y^0)],\quad \forall X^0, Y^0\in\Gamma(A_0).
 \end{eqnarray}
\end{lem}

When $n=1$, a Lie 1-algebroid is exactly a Lie algebroid. Associated to any vector bundle $E$, the covariant differential operator bundle $\dev (E)$ is a Lie algebroid naturally. Let $(A;a,l_2)$ be a Lie algebroid and $\nabla:A\longrightarrow \dev(E)$ a bundle map.  Then there is a differential operator $d_\nabla:\Gamma(\Hom(\wedge^kA,E))\longrightarrow \Gamma(\Hom(\wedge^{k+1}A,E))$ defined by
\begin{eqnarray*}
  d_{\nabla}\theta (X_1,\cdots, X_{k+1})&=&\sum_{i=1}^{k+1}(-1)^{i+1}\nabla_{X_i}\theta(X_1,\cdots,\widehat{X_i},\cdots, X_{k+1})\\
  &&+\sum_{i<j}(-1)^{i+j}\theta(l_2(X_i,X_j),X_1,\cdots,\widehat{X_i},\cdots,\widehat{X_j},\cdots, X_{k+1}),
\end{eqnarray*}
for all $X_1,\cdots, X_{k+1}\in\Gamma(A)$.  $d^2_\nabla=0$ if and only if $\nabla$ is a Lie algebroid morphism, i.e. a representation of $A$ on $E$. See \cite{General theory of Lie groupoid and Lie algebroid} for more details.

The notion of a $\LWX$ $2$-algebroid was introduced in \cite{LiuSheng} as the categorification of a Courant algebroid \cite{lwx,Roytenbergphdthesis}.

 \begin{defi}\label{defi:Courant-2 algebroid}
A {\bf $\LWX$ $2$-algebroid} is a graded vector bundle $\huaE=E_{-1}\oplus E_0$ over $M$ equipped with a nondegenerate symmetric bilinear form $S$ on $\huaE$, a bilinear operation $\diamond:\Gamma(E_{-i})\times \Gamma(E_{-j})\longrightarrow \Gamma(E_{-(i+j)})$, $0\leq i+j\leq 1$, which is skewsymmetric restricted on $\Gamma(E_0)\times \Gamma(E_0)$, an $E_{-1}$-valued $3$-form $\Omega$ on $E_0$, two bundle maps $\partial:E_{-1}\longrightarrow E_0$ and $\rho:E_0\longrightarrow TM$, such that $E_{-1}$ and $E_0$ are isotropic and the following conditions are satisfied:
\begin{itemize}
\item[$\rm(i)$]$(\Gamma(E_{-1}),\Gamma(E_0),\partial,\diamond,\Omega)$ is a Leibniz $2$-algebra,
\item[$\rm(ii)$]for all $e\in\Gamma(\huaE)$, $e\diamond e=\frac{1}{2}\huaD S(e,e)$,
\item[$\rm(iii)$]for all $e^1_1,e^1_2\in\Gamma(E_{-1})$, $S( \partial(e^1_1),e^1_2)=S(e^1_1,\partial(e^1_2))$,
\item[$\rm(iv)$]for all $e_1,e_2,e_3\in\Gamma(\huaE)$, $\rho(e_1)S( e_2,e_3)=S( e_1\diamond e_2,e_3)+S(e_2,e_1\diamond e_3)$,
\item[$\rm(v)$]for all $e^0_1,e^0_2,e^0_3,e^0_4\in\Gamma(E_0)$, $S(\Omega(e^0_1,e^0_2,e^0_3),e^0_4)=-S(e^0_3,\Omega(e^0_1,e^0_2,e^0_4))$,
 \end{itemize}
 where $\huaD:\CWM\longrightarrow \Gamma(E_{-1})$ is defined by
 \begin{equation}
 S(\huaD f,e^0)=\rho(e^0)(f),\quad \forall f\in\CWM, e^0\in\Gamma(E_0).
 \end{equation}

\end{defi}
Denote a $\LWX$ 2-algebroid  by $(E_{-1},E_0,\partial,\rho,S,\diamond,\Omega)$, or simply by $\huaE$. The skew-symmetrization of the Leibniz 2-algebra structure $\diamond$ will give rise to a Lie 3-algebra structure. There is a close relationship between $\LWX$ $2$-algebroids and symplectic NQ manifolds (QP manifolds) of degree 3. On the double of a Lie 2-bialgebroid, there is a $\LWX$ 2-algebroid structure naturally. See \cite{Ikeda,LiuSheng} for more details.

\section{Transitive Lie 2-algebroids}

In this section, we study the structure of a transitive Lie 2-algebroid. In Subsection 3.1, by choosing a splitting, we obtain some structure maps and we give the conditions that these structure maps satisfy. Then in Subsection 3.2, we characterize these conditions by a morphism from the tangent Lie algebroid $TM$ to a strict Lie 3-algebroid constructed from derivations of a graded bundle of Lie 2-algebras.

\subsection{General description of a transitive Lie 2-algebroid}\label{subsec:tranLie2}

 A Lie 2-algebroid $(A_{-1}\oplus A_0;\rho,l_1,l_2,l_3)$ is said to be {\bf transitive} if the anchor
  $\rho:A_0\longrightarrow TM$ is surjective.
Denote by $\huaG=\ker(\rho)$. By Lemma \ref{lem:propertyLien}, we deduce that $l_2(u,v)\in\Gamma(\huaG)$ for all $u,v\in\Gamma(\huaG)$. Then it is obvious that $(A_{-1}\oplus\huaG;\frkl_1,\frkl_2,\frkl_3)$ is a graded bundle of Lie 2-algebras, where $\frkl_1=l_1$, $\frkl_2$ and $\frkl_3$ are restrictions of $l_2$ and $l_3$ respectively.

\begin{defi}
  A splitting (connection) of a transitive Lie $2$-algebroid $(A_{-1}\oplus A_0;\rho,l_1,l_2,l_3)$ consists of  a section $\sigma:TM\longrightarrow A_0$ of the following short exact sequence of vector bundles:
\begin{equation}
  0\stackrel{}{\longrightarrow}\huaG\stackrel{\ii}{\longrightarrow}A_0\stackrel{\rho}{\longrightarrow}TM\stackrel{}{\longrightarrow}0,
\end{equation}
and a bundle map $\gamma:\wedge^2TM\longrightarrow A_{-1}$. Here $\ii$ denotes the inclusion map.
\end{defi}

After choosing a splitting, we have $A_0\cong TM\oplus\huaG$ and $\rho$ is simply the projection $\pr_{TM}$. Define $\nabla^0:TM\longrightarrow \dev(\huaG)$ and $\nabla^1:TM\longrightarrow \dev(A_{-1})$ by
\begin{equation}\label{eq:definabla}
  \nabla^0_Xu=l_2(\sigma(X),u),\quad \nabla^1_Xm=l_2(\sigma(X),m),\quad\forall X\in\frkX(M), u\in\Gamma(\huaG), m\in\Gamma(A_{-1}).
\end{equation}
Define the   bundle map $R:\wedge^2TM\longrightarrow \huaG$ by
\begin{equation}\label{eq:defiR}
  R(X,Y)=\sigma[X,Y]-l_2(\sigma(X),\sigma(Y))-\frkl_1\gamma(X,Y), \quad\forall X,Y\in\frkX(M).
\end{equation}
Define the   bundle map $I:\wedge^3TM\longrightarrow A_{-1}$ by
\begin{equation}\label{defi:I}
 I(X,Y,Z)=-d_{\nabla^1}\gamma(X,Y,Z)-l_3(\sigma(X),\sigma(Y),\sigma(Z)), \quad\forall X,Y,Z\in\frkX(M).
\end{equation}
Define totally skewsymmetric bundle map $J:\wedge^2TM\otimes\huaG\longrightarrow A_{-1}$ by
\begin{equation}\label{defi:J}
J(X,Y,u)=-l_3(\sigma(X),\sigma(Y),u), \quad\forall X,Y\in\frkX(M),\quad u\in\Gamma(\huaG).
\end{equation}
Define totally skewsymmetric bundle map $K:TM\otimes\wedge^2\huaG\longrightarrow A_{-1}$ by
\begin{equation}\label{defi:K}
K(X,u,v)=l_3(\sigma(X),u,v), \quad\forall X\in\frkX(M),\quad u,v\in\Gamma(\huaG).
\end{equation}

Transfer the transitive Lie 2-algebroid structure  from $A_{-1}\oplus A_0$ to $A_{-1}\oplus(\huaG\oplus TM)$, for which we use the same notation $(\rho,l_1,l_2,l_3)$, we have
\begin{equation}\label{eq:strformula}
\left\{\begin{array}{rcl}
\rho&=&\pr_{TM},\\
  l_1(m)&=&\frkl_1(m),\\
  l_2(X+u,Y+v)&=&[X,Y]-R_\gamma(X,Y)+\nabla^0_Xv-\nabla^0_Yu+\frkl_2(u,v),\\
  l_2(X+u,m)&=&-l_2(m,X+u)=\nabla^1_Xm+\frkl_2(u,m),\\
  l_3(X+u,Y+v,Z+w)&=&-I_\gamma(X,Y,Z)-J(X,Y,w)-J(X,v,Z)-J(u,Y,Z)\\
  &&+K(X,v,w)+K(u,Y,w)+K(u,v,Z)+\frkl_3(u,v,w),
  \end{array}\right.
\end{equation}
where $$R_\gamma=R+\frkl_1\circ \gamma,\quad I_\gamma=I+d_{\nabla^1}\gamma.$$

\begin{rmk}
 $R_\gamma$ and $I_\gamma$ can be viewed as curvatures.
\end{rmk}

\begin{thm}\label{thm:transtiveLie2}
Let $(A_{-1}\oplus\huaG;\frkl_1,\frkl_2,\frkl_3)$ be a graded bundle of Lie $2$-algebras. Then  $(A_{-1}\oplus(\huaG\oplus TM);\rho,l_1,l_2,l_3)$ is a transitive Lie $2$-algebroid, where $\rho,l_1,l_2,l_3$ are given by \eqref{eq:strformula} for totally skew-symmetric $I,J,K$, if and only if for all $X,Y,Z\in\frkX(M)$, $u,v,w\in\Gamma(\huaG)$ and $m\in\Gamma(A_{-1})$, the following equalities hold:
  \begin{eqnarray}
   \label{eq:compatl1} \frkl_1\circ \nabla^1_X&=&\nabla^0_X\circ \frkl_1,\\
 \label{eq:thmder1} \nabla^0_X\frkl_2(v,w)-\frkl_2(\nabla^0_Xv,w)-\frkl_2(v,\nabla^0_Xw)&=&\frkl_1 K(X,v,w),\\
     \label{eq:thmder2} \nabla^1_X\frkl_2(v,m)-\frkl_2(\nabla^0_Xv,m)-\frkl_2(v,\nabla^1_Xm)&=& K(X,v,\frkl_1 m),\\
    \nonumber \frkl_2(u,K(v,w,X))-\frkl_2(v,K(u,w,X))+\frkl_2(w,K(u,v,X))&&\\
     \nonumber -K(\frkl_2(u,v),w,X)-K(v,\frkl_2(u,w),X)+K(u,\frkl_2(v,w),X)&&\\
     \label{eq:thmder3}  -\nabla^1_X\frkl_3(u,v,w)+\frkl_3(\nabla^0_Xu, v,w)+\frkl_3(u,\nabla^0_Xv,w)+\frkl_3(u,v,\nabla^0_Xw)&=&0,\\
       \label{eq:thmobsmor1} \nabla_X^0\nabla^0_Yw- \nabla_Y^0\nabla^0_Xw- \nabla_{[X,Y]}^0 w+\frkl_2(R_\gamma(X,Y),w)+\frkl_1 J(X,Y,w)&=&0,\\
   \label{eq:thmobsmor2} \nabla_X^1\nabla^1_Ym- \nabla_Y^1\nabla^1_Xm- \nabla_{[X,Y]}^1m+\frkl_2(R_\gamma(X,Y),m)+J(X,Y,\frkl_1m)&=&0,\\
        \label{eq:thmRcon}d_{\nabla^0}R_\gamma(X,Y,Z)&=&\frkl_1 I_\gamma(X,Y,Z),\\
             \nonumber \frkl_2(u,I_\gamma(X,Y,Z))-\nabla^1_XJ(u,Y,Z)+\nabla^1_YJ(u,X,Z)&&\\-\nabla^1_ZJ(u,X,Y)
      \nonumber+J(\nabla^0_Xu,Y,Z) +J(X,\nabla^0_Yu, Z)+J(X,Y,\nabla^0_Zu)&&\\ \nonumber+J(u,[X,Y],Z)+J(u,[Z,X],Y)+J(u,[Y,Z],X)&&\\
     \label{eq:thmJcon}-K(u,Z,R_\gamma(X,Y))-K(u,Y,R_\gamma(Z,X))-K(u,X,R_\gamma(Y,Z))&=&0,\\
     \label{eq:thmIcon}d_{\nabla^1}I_\gamma+J\circ R_\gamma&=&0,
  \end{eqnarray}
 where  $J\circ R_\gamma:\wedge^4\frkX(M)\longrightarrow \Gamma(A_{-1})$ is given by
  \begin{equation}
    J\circ R_\gamma(X_1,\cdots, X_4)=\frac{1}{4}\sum_{\tau\in S_4} \sgn(\tau)J(R_\gamma(X_{\tau(1)},X_{\tau(2)}),X_{\tau(3)} ,X_{\tau(4)}).
  \end{equation}
\end{thm}

\pf Assume that $(A_{-1}\oplus(\huaG\oplus TM);\rho,l_1,l_2,l_3)$ is a transitive Lie $2$-algebroid.   By the equality $l_1l_2(X,m)=l_2(X,l_1(m))$, we deduce that \eqref{eq:compatl1} holds. By the equality
$$l_2(X,l_2(v,w))+l_2(v,l_2(w,X))+l_2(w,l_2(X,v))=l_1l_3(X,v,w),$$
we deduce that \eqref{eq:thmder1} holds. Similarly, we deduce that \eqref{eq:thmder2} holds. By the equality
{\footnotesize\begin{eqnarray*}
 l_2(u,l_3(v,w,X))-l_2(v,l_3(u,w,X))+l_2(w,l_3(u,v,X))-l_2(X,l_3(u,v,w))-l_3(l_2(u,v),w,X)\\
 + l_3(l_2(u,w),v,X)-l_3(l_2(v,w),u,X)-l_3(l_2(u,X),v,w)+l_3(l_2(v,X),u,w)-l_3(l_2(w,X),u,v)=0,
\end{eqnarray*}
}
we deduce that \eqref{eq:thmder3} holds. By the equality
$$
l_2(X,l_2(Y,w))+l_2(Y,l_2(w,X))+l_2(w,l_2(X,Y))=l_1l_3(X,Y,w),
$$
we deduce that \eqref{eq:thmobsmor1} holds. Similarly, we deduce that \eqref{eq:thmobsmor2} holds. By the equality
$$
l_2(X,l_2(Y,Z))+l_2(Y,l_2(Z,X))+l_2(Z,l_2(X,Y))=l_1l_3(X,Y,Z),
$$
we deduce that \eqref{eq:thmRcon} holds. Then by the equality
{\footnotesize\begin{eqnarray*}
l_2(u,l_3(X,Y,Z))-l_2(X,l_3(u,Y,Z))+l_2(Y,l_3(u,X,Z))-l_2(Z,l_3(u,X,Y))-l_3(l_2(u,X),Y,Z) \\
+l_3(l_2(u,Y),X,Z)-l_3(l_2(u,Z),X,Y)-l_3(l_2(X,Y),u,Z)+l_3(l_2(X,Z),u,Y)-l_3(l_2(Y,Z),u,X)=0,
\end{eqnarray*}
}
we deduce that \eqref{eq:thmJcon} holds. Finally, by the equality
{\footnotesize
\begin{eqnarray*}
 \sum_{i=1}^{4}(-1)^{i+1}l_2({X_i},l_3(X_1,\cdots,\widehat{X_i},\cdots, X_{4}))+\sum_{i<j}(-1)^{i+j}l_3(l_2(X_i,X_j),X_1,\cdots,\widehat{X_i},\cdots,\widehat{X_j},\cdots, X_{4})=0,
\end{eqnarray*}
}
we deduce that \eqref{eq:thmIcon} holds.

Conversely, by \eqref{eq:thmder1} and \eqref{eq:thmobsmor1}, we can deduce that
  \begin{eqnarray}
      \nonumber -\frkl_2(u,J(v,Y,Z))+\frkl_2(v,J(u,Y,Z))+J(\frkl_2(u,v),Y,Z)&&\\
      \nonumber+\nabla^1_YK(u,v,Z)-\nabla^1_ZK(u,v,Y)-K(u,v,[Y,Z])+\frkl_3(u,v,R_\gamma(Y,Z))&&\\
    \label{eq:thmKcon} +K(v,\nabla^0_Yu,Z)+K(v,Y,\nabla_Z^0u)-K(u,\nabla_Y^0v,Z)-K(u,Y,\nabla^0_Zv)&=&0.
     \end{eqnarray}
 Then by \eqref{eq:compatl1}-\eqref{eq:thmIcon} and \eqref{eq:thmKcon}, we can deduce that  $(A_{-1}\oplus (\huaG\oplus TM);\rho,l_1,l_2,l_3)$ is a transitive Lie $2$-algebroid. \qed

\subsection{Transitive Lie 2-algebroids and Lie 3-algebroid morphisms}

In this subsection, we give a conceptual explanation of equations listed in Theorem \ref{thm:transtiveLie2}. We show that they give rise to a morphism from the Lie algebroid $TM$ to a strict Lie 3-algebroid constructed from derivations of the graded bundle of Lie 2-algebras $(A_{-1}\oplus\huaG;\frkl_1,\frkl_2,\frkl_3)$.

First we give the definition of a morphism from a Lie algebroid to a strict Lie 3-algebroid. See \cite[Section 4.1]{BP} for details for general formulas of a morphism between Lie $n$-algebroids.

\begin{defi}\label{Def2}
Let $\huaA=(A;\rho,l_2)$ be a Lie algebroid and $\huaA'=(A_{-2}'\oplus A_{-1}'\oplus A_0';\rho',l_1',l_{2}')$
a strict Lie $3$-algebroid. A morphism $F$ from $\huaA$ to $\huaA'$ consists of:
\begin{itemize}
\item[$\bullet$] a bundle map $F^{1}:A\longrightarrow A_{0}'$,
\item[$\bullet$] a bundle map $F^{2}:\wedge^2 A_{0} \longrightarrow  A_{-1}'$,
 \item[$\bullet$] a bundle map $F^{3}: \wedge^3  A_0 \longrightarrow  A_{-2}',$
\end{itemize}
such that for all $ X,Y,Z,X_i\in \Gamma(A)$, $i=1,2,3,4$,  we have
\begin{eqnarray*}
\rho'\circ F^1&=&\rho,\\
F^{1}l_{2}(X,Y)-l_{2}'(F^{1}(X),F^{1}(Y))&=&l_1'F^{2}(X,Y),\\
l_{2}'(F^{1}(X),F^{2}(Y,Z))-F^{2} (l_{2}(X,Y),Z)+c.p.&=&l_1'F^{3}(X,Y,Z),
\end{eqnarray*}
 and
\begin{eqnarray*}
&&\sum_{i=1}^4(-1)^{i+1}l_2'(F^1(X_i),F^3(X_1,\cdots,\widehat{X_i},\cdots X_4))\\
&&+\sum_{i<j}(-1)^{i+j}\Big(F^3(l_2(X_i,X_j),X_k,X_l)-\half l_2'(F^2(X_i,X_j),F^2(X_k,X_l))\Big)=0,
\end{eqnarray*}
where $k<l$ and $\{k,l\}\cap \{i,j\}=\emptyset.$
\end{defi}

Then we construct a strict Lie 3-algebroid using derivations of a graded bundle of Lie 2-algebras $\frkG=(\huaG_{-1}\oplus \huaG_0;\frkl_1,\frkl_2,\frkl_3)$. Denote by $\dev(\huaG_0)$ and $\dev(\huaG_{-1})$ the covariant differential operator bundle associated to
$\huaG_0$ and $\huaG_{-1}$ respectively.  We have the following exact sequences:
\begin{eqnarray*}
\xymatrix@C=0.5cm{0 \ar[r] & \gl(\huaG_0)  \ar[rr] && \dev(\huaG_0)  \ar[rr]^{\jd_0} && TM \ar[r]  & 0,}\\
\xymatrix@C=0.5cm{0 \ar[r] & \gl(\huaG_{-1})  \ar[rr]  && \dev(\huaG_{-1})  \ar[rr]^{\jd_{1}} && TM \ar[r]  & 0.}
\end{eqnarray*}
\begin{defi}
  A {\bf degree $0$ derivation} of a graded bundle of Lie $2$-algebras $(\huaG_{-1}\oplus\huaG_0;\frkl_1,\frkl_2,\frkl_3)$ is a triple $\frkd=(\frkd_0,\frkd_1,l_\frkd)$, where $\frkd_0\in\Gamma(\dev (\huaG_0))$, $\frkd_1\in\Gamma(\dev (\huaG_{-1}))$ and $l_\frkd\in\Gamma(\Hom(\wedge^2\huaG_0, \huaG_{-1}))$, such that we have
 \begin{eqnarray*}
 \jd_0(\frkd_0)&=&\jd_{1}(\frkd_1),\\
   \frkd_0\circ \frkl_1&=&\frkl_1\circ \frkd_1,\\
   \frkd_0\frkl_2(u,v)&=&\frkl_2(\frkd_0(u),v)+\frkl_2(u,\frkd_0(v))+\frkl_1 l_\frkd(u,v),\\
   \frkd_1\frkl_2(u,m)&=&\frkl_2(\frkd_0(u),m)+\frkl_2(u,\frkd_1(m))+ l_\frkd(u,\frkl_1 (m)),
   \end{eqnarray*}
   and
    \begin{eqnarray*}
   &&l_{\frkd}(u,\frkl_2(v,w))+\frkl_2(u,l_{\frkd}(v,w))+\frkl_{3} (\frkd_0(u),v,w)+\frkl_{3} (u,\frkd_0(v),w)+\frkl_{3} (u,v,\frkd_0(w)) \nonumber \\
&=&\frkd_1 \frkl_{3}(u,v,w)+l_{\frkd}(\frkl_2(u,v),w)+l_{\frkd}(v,\frkl_2(u,w))+\frkl_2(l_{\frkd}(u,v),w)+\frkl_2(v,l_{\frkd}(u,w)),
 \end{eqnarray*}
 for all $u,v,w\in\Gamma(\huaG_0)$ and $m\in\Gamma(\huaG_{-1})$.
\end{defi}
It is obvious that the set of degree 0 derivations is the section space of a vector bundle, which we denote by $\Der^0(\frkG)$. Furthermore, we define $\Der^{-1}(\frkG)=\Hom(\huaG_0,\huaG_{-1})$.

\begin{ex}
  For any $u\in\Gamma(\huaG_0)$, define $\ad^0_u:\Gamma(\huaG_0)\longrightarrow \Gamma(\huaG_0)$ and $\ad^{1}_u:\Gamma(\huaG_{-1})\longrightarrow \Gamma(\huaG_{-1})$ by
  $$
  \ad^0_uv=\frkl_2(u,v),\quad\ad^{1}_um=\frkl_2(u,m),\quad \forall v\in\Gamma(\huaG_{0}), m\in\Gamma(\huaG_{-1}).
  $$
  Define $l_{\ad_u}:\wedge^2\Gamma(\huaG_0)\longrightarrow\Gamma(\huaG_{-1})$ by
  $$
  l_{\ad_u}(v,w)=\frkl_3(u,v,w),\quad \forall v,w\in\Gamma(\huaG_{0}).
  $$
  Then it is straightforward to see that $(\ad^0_u,\ad^{1}_u,l_{\ad_u})\in\Gamma(\Der^0(\frkG))$.
\end{ex}

 Define $\jd:\Der^0(\frkG)\longrightarrow TM$  by
$$
\jd(\frkd_0,\frkd_1,l_{\frkd}):=\jd_0(\frkd_0).
$$
 On the graded bundle $\Der(\frkG)=\Der^0(\frkG)\oplus \Der^{-1}(\frkG)$, we define $\delta:\Der^{-1}(\frkG)\longrightarrow \Der^0(\frkG)$ and  skew-symmetric bracket operation $[\cdot,\cdot]_{\Der}:\Gamma(\Der^{-i}(\frkG))\times \Gamma(\Der^{-j}(\frkG))\longrightarrow \Gamma(\Der^{-(i+j)}(\frkG))$, $0\leq i+j\leq 1$, by
\begin{eqnarray}
  \delta(\phi)&=&(\frkl_1\circ \phi,\phi\circ \frkl_1,l_{\delta(\phi)}),\\
  ~[(\frkd_0,\frkd_1,l_\frkd),(\frkt_0,\frkt_1,l_\frkt)]_\Der&=&([\frkd_0,\frkt_0]_C,[\frkd_1,\frkt_1]_C,l_{[\frkd,\frkt]_C}),\\
  ~ [(\frkd_0,\frkd_1,l_\frkd),\phi]_\Der&=&\frkd_1\circ \phi-\phi\circ\frkd_0,
\end{eqnarray}
where  $l_{\delta(\phi)}, l_{[\frkd,\frkt]_C}\in\Gamma(\Hom(\wedge^2\huaG_0,\huaG_{-1}))$ are defined by
\begin{eqnarray*}
l_{\delta(\phi)}(u,v)&=&\phi(\frkl_2(u,v))-\frkl_2(\phi(u),v)-\frkl_2(u,\phi(v)),\\
  l_{[\frkd,\frkt]_C}(u,v)&=&\frkd_1 l_\frkt(u,v)-l_\frkt(\frkd_0(u),v)-l_\frkt(u,\frkd_0(v))    + l_\frkd(\frkt_0(u),v)+l_\frkd(u,\frkt_0(v))-\frkt_1 l_\frkd(u,v).
\end{eqnarray*}

It is straightforward to deduce that
\begin{pro}
 Let  $\frkG=(\huaG_{-1}\oplus\huaG_0;\frkl_1,\frkl_2,\frkl_3)$ be a graded bundle of Lie $2$-algebras. Then $(\Der^{-1}(\frkG)\oplus\Der^0(\frkG);\jd,\delta,[\cdot,\cdot]_\Der)$ is a strict Lie $2$-algebroid.
\end{pro}

We go on constructing a strict Lie 3-algebroid. On the 3-term graded vector bundles
$$
\DER(\frkG)=\huaG_{-1}\oplus(\Der^{-1}(\frkG)\oplus \huaG_0)\oplus \Der^0(\frkG),
$$
where the degree 0 part $\DER^0(\frkG)$ is $\Der^0(\frkG)$, the degree $-1$ part $\DER^{-1}(\frkG)$ is $\Der^{-1}(\frkG)\oplus \huaG_0$ and the degree $-2$ part $\DER^{-2}(\frkG)$ is $\huaG_{-1}$, we define $\dM:\DER^{-i}(\frkG)\longrightarrow \DER^{-i+1}(\frkG)$, $i=1,2$ and $[\cdot,\cdot]_\DER:\Gamma(\DER^{-i}(\frkG))\times \Gamma(\DER^{-j}(\frkG))\longrightarrow \Gamma(\DER^{-(i+j)}(\frkG))$, $0\leq i+j\leq 2$, by
\begin{eqnarray*}
  \dM(m)&=&-\frkl_2(m,\cdot)+\frkl_1(m),\\
  \dM(\phi+u)&=&\delta(\phi)+(\ad^0_u,\ad^{1}_u,l_{\ad_u}),\\
 ~ [(\frkd_0,\frkd_1,l_\frkd),(\frkt_0,\frkt_1,l_\frkt)]_\DER&=& [(\frkd_0,\frkd_1,l_\frkd),(\frkt_0,\frkt_1,l_\frkt)]_\Der,\\
 ~  [(\frkd_0,\frkd_1,l_\frkd),\phi+u]_\DER&=&[(\frkd_0,\frkd_1,l_\frkd),\phi]_\Der+l_\frkd(u,\cdot)+\frkd_0(u),\\
 ~[(\frkd_0,\frkd_1,l_\frkd),m]_\DER&=&\frkd_1(m),\\
 ~[\phi+u,\psi+v]_\DER&=&\phi(v)+\psi(u),
\end{eqnarray*}
for all $(\frkd_0,\frkd_1,l_\frkd),(\frkt_0,\frkt_1,l_\frkt)\in\Gamma(\Der^0(\frkG))$, $\phi,\psi\in\Gamma(\Der^{-1}(\frkG))$, $u,v\in\Gamma(\huaG_0)$ and $m\in\Gamma(\huaG_{-1})$.

\begin{thm}
   Let  $\frkG=(\huaG_{-1}\oplus\huaG_0;\frkl_1,\frkl_2,\frkl_3)$ be a graded bundle of Lie $2$-algebras. Then $$(\DER^{-2}(\frkG)\oplus\DER^{-1}(\frkG)\oplus\DER^0(\frkG);\jd,\dM,[\cdot,\cdot]_\DER)$$ is a strict Lie $3$-algebroid.
\end{thm}
\pf The proof is a straightforward verification. We leave it to readers. \qed \vspace{3mm}

With above preparations, we go back to Subsection 3.1. For a transitive Lie 2-algebroid $(A_{-1}\oplus A_0;\rho,l_1,l_2,l_3)$, let $\frkG=(A_{-1}\oplus\huaG;\frkl_1,\frkl_2,\frkl_3)$ be the corresponding graded bundle of Lie 2-algebras, where $\huaG=\ker(\rho)$, and $\DER(\frkG)$ the strict Lie 3-algebroid given above.

Define $F^1:TM\longrightarrow \DER^0(\frkG)$ by
\begin{equation}
  F^1(X)=(\nabla^0_X,\nabla^1_X,K(X,\cdot,\cdot)), \quad\forall X\in\frkX(M).
\end{equation}
By \eqref{eq:compatl1}, \eqref{eq:thmder1}, \eqref{eq:thmder2}, \eqref{eq:thmder3}, $F^1$ is well-defined.

Define $F^2:\wedge^2TM\longrightarrow \DER^{-1}(\frkG)$ by
\begin{equation}
  F^2(X,Y)=J(X,Y,\cdot)+R_\gamma(X,Y), \quad\forall X,Y\in\frkX(M).
\end{equation}
Define $F^3:\wedge^3TM\longrightarrow \DER^{-2}(\frkG)$ by
\begin{equation}
  F^3(X,Y,Z)=I_\gamma(X,Y,Z), \quad\forall X,Y,Z\in\frkX(M).
\end{equation}

\begin{thm}
Let $(A_{-1}\oplus A_0;\rho,l_1,l_2,l_3)$ be a transitive Lie $2$-algebroid.  Then $(F^1,F^2,F^3)$ is a morphism from the Lie algebroid $TM$ to the strict Lie $3$-algebroid $\DER(\frkG)$.
\end{thm}
\pf We only give a sketch of the proof and leave details to readers. By
\eqref{eq:thmobsmor1}, \eqref{eq:thmobsmor2} and \eqref{eq:thmKcon}, we deduce that
$$
F^1([X,Y])-[F^1(X),F^1(Y)]_\DER=\dM F^2(X,Y),\quad\forall X,Y\in\frkX(M).
$$
By \eqref{eq:thmRcon} and \eqref{eq:thmJcon}, we deduce that
$$
[F^{1}(X),F^{2}(Y,Z)]_\DER-F^{2} ([X,Y],Z)+c.p.=\dM F^{3}(X,Y,Z).
$$
By \eqref{eq:thmIcon}, we deduce that
\begin{eqnarray*}
&&\sum_{i=1}^4(-1)^{i+1}[F^1(X_i),F^3(X_1,\cdots,\widehat{X_i},\cdots X_4)]_\DER\\
&&+\sum_{i<j}(-1)^{i+j}\Big(F^3([X_i,X_j],X_k,X_l)-\half [F^2(X_i,X_j),F^2(X_k,X_l)]_\DER\Big)\\
&=&(d_{\nabla^1}I_\gamma+J\circ R_\gamma)(X_1,\cdots, X_4)\\
&=&0.
\end{eqnarray*}
Thus, $(F^1,F^2,F^3)$ is a morphism form $TM$ to $\DER(\frkG)$.\qed

\section{The first Pontryagin class of a quadratic Lie 2-algebroid }

In this section, we give the notion of a quadratic Lie 2-algebroid and define its first Pontryagin class, which is a cohomology class in $H^5(M)$.

\begin{defi}\label{defi:quadraticLie2}
  A transitive Lie $2$-algebroid $(A_{-1}\oplus A_0;\rho,l_1,l_2,l_3)$ is said to be a {\bf quadratic Lie $2$-algebroid} if there is a degree $1$ nondegenerate graded symmetric bilinear form $\huaS$ on the graded vector bundle $A_{-1}\oplus \ker(\rho)$ satisfying
  \begin{eqnarray}
    \label{eq:Lie2inv1}\huaS(\frkl_1(m),p)&=&\huaS(m,\frkl_1(p)),\\
    \label{eq:Lie2inv2}\rho(e^0)\huaS(u,m)&=&\huaS(l_2(e^0,u),m)+\huaS(u,l_2(e^0,m)),\\
    \label{eq:Lie2inv3}\huaS(l_3(e^0_1,e^0_2,u),v)&=&-\huaS(l_3(e^0_1,e^0_2,v),u),
  \end{eqnarray}
  for all $u,v\in\Gamma(\ker(\rho))$, $e^0,e^0_1,e^0_2\in\Gamma(A_0)$ and $m,p\in\Gamma(A_{-1})$.
\end{defi}

Let $(A_{-1}\oplus A_0;\rho,l_1,l_2,l_3)$  be a   quadratic Lie $2$-algebroid. As before, denote by $\huaG=\ker(\rho)$. By the nondegeneracy of $\huaS$, we deduce that $A_{-1}\cong \huaG^*$ and we will write $A_{-1}=\huaG^*$ directly. Then the pairing $\huaS$ is simply given by
\begin{eqnarray}\label{defi:huaS}
  \huaS(u+m,v+p)=\langle u,p\rangle+\langle v,m\rangle,\quad \forall u,v\in\Gamma(\huaG), m,p\in\Gamma(\huaG^*).
\end{eqnarray}
By \eqref{eq:Lie2inv1}-\eqref{eq:Lie2inv3}, we deduce that for the graded bundle of Lie 2-algebras $(\huaG^*\oplus\huaG;\frkl_1,\frkl_2,\frkl_3)$, there holds:
\begin{eqnarray}
 \label{eq:Lie2inv11} \langle \frkl_1(m),p\rangle&=&\langle m,\frkl_1(p)\rangle, \quad \forall m,p\in\Gamma(\huaG^*),\\
 \label{eq:Lie2inv22} \langle\frkl_2(u,v),m\rangle&=&-\langle v,\frkl_2(u,m)\rangle,\quad \forall u,v\in\Gamma(\huaG), m\in\Gamma(\huaG^*),\\
 \label{eq:Lie2inv33}\langle \frkl_3(u,v,w),x \rangle&=&-\langle w,\frkl_3(u,v,x) \rangle,\quad \forall u,v,w,x\in\Gamma(\huaG).
\end{eqnarray}

\begin{pro}\label{pro:quadraticLie2}
With the same assumption and conditions in Theorem \ref{thm:transtiveLie2}, the transitive Lie $2$-algebroid $(\huaG^*\oplus (\huaG\oplus TM);\rho,l_1,l_2,l_3)$ is a quadratic Lie $2$-algebroid, where $\rho,l_1,l_2,l_3$ are given by \eqref{eq:strformula} for totally skew-symmetric $I,J,K$, if and only if \eqref{eq:Lie2inv11}-\eqref{eq:Lie2inv33} and the following equalities hold:
  \begin{eqnarray}
   \label{eq:thminv1}\langle\nabla^0_Xu,m\rangle+\langle u,\nabla^1_Xm\rangle&=&X\langle u,m\rangle,\\
   \label{eq:thminv2} \langle J(X,Y,u),v\rangle+\langle J(X,Y,v),u\rangle&=&0,\\
    \label{eq:thminv3}\langle K(X,u,v),w\rangle+\langle K(X,u,w),v\rangle&=&0.
\end{eqnarray}
\end{pro}
\pf It follows from that the invariant conditions in the definition of a quadratic Lie $2$-algebroid is equivalent to \eqref{eq:Lie2inv11}-\eqref{eq:thminv3}. \qed\vspace{3mm}

Let $(\huaG^*\oplus(\huaG\oplus TM);\rho,l_1,l_2,l_3)$ be a quadratic Lie $2$-algebroid given by \eqref{eq:strformula} with $\huaS$ the graded symmetric bilinear form given by \eqref{defi:huaS}. Define $\huaS(R_\gamma,I_\gamma)\in\Omega^5(M)$  by
  \begin{equation}\label{defi:5form}
    \huaS(R_\gamma,I_\gamma)(X_1,\cdots, X_5)=\frac{1}{12}\sum_{\tau\in S_5}\sgn(\tau)\huaS\big(R_\gamma(X_{\tau(1)},X_{\tau(2)}),I_\gamma(X_{\tau(3)},X_{\tau(4)},X_{\tau(5)})\big).
  \end{equation}

  \begin{thm}\label{thm:5form}
    The $5$-form $\huaS(R_\gamma,I_\gamma)$ is closed, i.e. $d\huaS(R_\gamma,I_\gamma)=0$, and its cohomology class $[\huaS(R_\gamma,I_\gamma)]$ in $H^5(M)$ does not depend on the choices of $\sigma$ and $\gamma$.
  \end{thm}

  The cohomology class $[\huaS(R_\gamma,I_\gamma)]\in H^5(M)$ is called {\bf the first Pontryagin class} of the quadratic Lie 2-algebroid $(A_{-1}\oplus A_0;\rho,l_1,l_2,l_3)$. To prove the above theorem, we need the following technical lemma.

  \begin{lem}\label{lem:d}
    For all  $\Xi\in\Gamma(\Hom(\wedge^kTM,\huaG))$ and  $\Pi\in\Gamma(\Hom(\wedge^lTM,\huaG^*))$, we have
    \begin{equation}\label{eq:dsep}
      d\huaS(\Xi,\Pi)=\huaS(d_{\nabla^0}\Xi,\Pi)+(-1)^k\huaS(\Xi,d_{\nabla^1}\Pi),
    \end{equation}
    where $\huaS(\Xi,\Pi)\in\Omega^{k+l}(M)$ is defined by
    $$
    \huaS(\Xi,\Pi)(X_1,\cdots, X_{k+l})=\frac{1}{k!l!}\sum_{\tau\in S_{k+l}}\sgn(\tau)\huaS\big(\Xi(X_{\tau(1)},\cdots,X_{\tau(k)}),\Pi(X_{\tau(k+1)},\cdots,X_{\tau(k+l)})\big),
$$
for all $X_1,\cdots, X_{k+l}\in\frkX(M)$.
  \end{lem}
\pf First we have
\begin{eqnarray*}
  && d\huaS(\Xi,\Pi)(X_1,\cdots, X_{k+l+1})\\&=&\sum_{i=1}^{k+l+1}(-1)^{i+1}X_i\huaS(\Xi,\Pi)(X_1,\cdots,\hat{X_i},\cdots,X_{k+l+1})\\
    &&+\sum_{i<j}(-1)^{i+j}\huaS(\Xi,\Pi)([X_i,X_j],X_1,\cdots,\hat{X_i},\cdots,\hat{X_j},\cdots,X_{k+l+1}).
    \end{eqnarray*}
    Then by \eqref{eq:thminv1}, we can deduce that \eqref{eq:dsep} holds after a careful computation. We omit details.\qed\vspace{3mm}

{\bf The proof of Theorem \ref{thm:5form}:} By Lemma \ref{lem:d}, \eqref{eq:thmRcon}, \eqref{eq:thmIcon}, \eqref{eq:Lie2inv11} and  \eqref{eq:thminv2}, we have
\begin{eqnarray*}
 d \huaS(R_\gamma,I_\gamma)&=&\huaS(d_{\nabla^0}R_\gamma,I_\gamma)+\huaS(R_\gamma,d_{\nabla^1}I_\gamma)\\
 &=&\huaS(\frkl_1 I_\gamma,I_\gamma)-\huaS(R_\gamma,J\circ R_\gamma)\\
 &=&0.
\end{eqnarray*}
Thus, $\huaS(R_\gamma,I_\gamma)$ is a closed 5-form.

Furthermore, by Lemma \ref{lem:d} and \eqref{eq:thmRcon}, we have
\begin{eqnarray*}
  \huaS(R_\gamma,I_\gamma)&=&\huaS(R,I)+\huaS(R,d_{\nabla^1}\gamma)+\huaS(\frkl_1\circ\gamma,I)+\huaS(\frkl_1\circ\gamma,d_{\nabla^1}\gamma)\\
  &=&\huaS(R,I)+d\huaS(R,\gamma)+\half d\huaS(\frkl_1\circ\gamma,\gamma),
\end{eqnarray*}
which implies that the cohomology class does not depend on the choices of $\gamma$.

If we choose another section $\sigma':TM\longrightarrow A_0$, then we define $\theta:TM\longrightarrow \huaG$ by
$$
\theta(X)=\sigma(X)-\sigma'(X),\quad\forall X\in\frkX(M).
$$
We have
\begin{eqnarray*}
R'(X,Y)&=&\sigma'[X,Y]-l_2(\sigma'(X),\sigma'(Y))-\frkl_1\circ \gamma(X,Y)\\
&=&R(X,Y)+d_{\nabla^0}\theta(X,Y)-\frkl_2(\theta(X),\theta(Y))
\end{eqnarray*}
and
\begin{eqnarray*}
I'(X,Y,Z)&=&-l_3(\sigma'(X),\sigma'(Y),\sigma'(Z))-d_{\nabla^1}\gamma(X,Y,Z)\\
&=&I(X,Y,Z)+\frkl_3(\theta(X),\theta(Y),\theta(Z))
-\Big(J(X,Y,\theta(Z))+K(X,\theta(Y),\theta(Z))+c.p.\Big).
\end{eqnarray*}
We define $\frkl_2\circ\theta:\wedge^2TM\longrightarrow\huaG$, $\frkl_3\circ\theta:\wedge^3TM\longrightarrow\huaG^*$, $J\circ\theta:\wedge^3TM\longrightarrow\huaG^*$ and $K\circ\theta:\wedge^3TM\longrightarrow\huaG^*$ respectively by
\begin{eqnarray*}
  \frkl_2\circ\theta(X,Y)&=&\frkl_2(\theta(X),\theta(Y)),\\
  \frkl_3\circ\theta(X,Y,Z)&=&\frkl_3(\theta(X),\theta(Y),\theta(Z)),\\
  J\circ\theta(X,Y,Z)&=&J(\theta(X),Y,Z)+J(X,\theta(Y),Z)+J(X,Y,\theta(Z)),\\
  K\circ\theta(X,Y,Z)&=&K(X,\theta(Y),\theta(Z))+K(\theta(X),Y,\theta(Z))+K(\theta(X),\theta(Y),Z).
\end{eqnarray*}
Therefore, by the following Lemma \ref{lem:thmI}-Lemma \ref{lem:thml2l3}, we have
\begin{eqnarray*}
\huaS(R'_\gamma,I'_\gamma)&=&\huaS(R_\gamma+d_{\nabla^0}\theta-\frkl_2\circ\theta,I_\gamma-J\circ\theta-K\circ\theta+\frkl_3\circ\theta)\\
&=&\huaS(R_\gamma,I_\gamma)+\huaS(d_{\nabla^0}\theta,I_\gamma)-\huaS(R_\gamma,J\circ\theta)-\huaS(R_\gamma,K\circ\theta)+\huaS(R_\gamma,\frkl_3\circ\theta)\\
&&-\huaS(d_{\nabla^0}\theta,J\circ\theta)-\huaS(d_{\nabla^0}\theta,K\circ\theta)+\huaS(d_{\nabla^0}\theta,\frkl_3\circ\theta)\\
&&-\huaS(\frkl_2\circ\theta,I_\gamma)+\huaS(\frkl_2\circ\theta,J\circ\theta)+\huaS(\frkl_2\circ\theta,K\circ\theta)-\huaS(\frkl_2\circ\theta,\frkl_3\circ\theta)\\
&=&\huaS(R_\gamma,I_\gamma)+d(C_4^I-C_4^J-C_4^K+C_4^{\frkl_3}),
\end{eqnarray*}
which implies that $\huaS(R'_\gamma,I'_\gamma)$ and $\huaS(R_\gamma,I_\gamma)$ are in the same cohomology class.
\qed\vspace{3mm}

For any $\huaG$-valued 1-form $\theta\in\Gamma(\Hom(TM,\huaG))$, we define a 4-form $C_4^{I}\in\Omega^4(M)$ by
\begin{equation}
  C_4^{I}= \huaS(\theta,I_\gamma).
\end{equation}
\begin{lem}\label{lem:thmI}
  With the above notations, we have
  \begin{equation}
    dC_4^I=\huaS(d_{\nabla^0}\theta,I_\gamma)-\huaS(R_\gamma,J\circ\theta).
  \end{equation}
\end{lem}
\pf By Lemma \ref{lem:d}, \eqref{eq:thmIcon} and \eqref{eq:thminv2}, we have
\begin{eqnarray*}
 d\huaS(\theta,I_\gamma)&=&\huaS(d_{\nabla^0}\theta,I_\gamma)-\huaS(\theta,d_{\nabla^1}I_\gamma) \\
 &=&\huaS(d_{\nabla^0}\theta,I_\gamma)+\huaS(\theta,J\circ R_\gamma)\\
 &=&\huaS(d_{\nabla^0}\theta,I_\gamma)-\huaS(R_\gamma,J\circ\theta).
\end{eqnarray*}
The proof is finished. \qed\vspace{3mm}

For any $\huaG$-valued 1-form $\theta\in\Gamma(\Hom(TM,\huaG))$, we define a 4-form $C_4^{J}\in\Omega^4(M)$ by
\begin{equation}
  C_4^{J}=\frac{1}{2}\huaS(\theta,J\circ\theta).
\end{equation}
\begin{lem}\label{lem:thmJ}
  With the above notations, we have
  \begin{equation}
    dC_4^J=\huaS(d_{\nabla^0}\theta,J\circ \theta)+\huaS(R_\gamma,K\circ\theta)+\huaS(\frkl_2\circ\theta,I_\gamma).
  \end{equation}
\end{lem}
\pf  By Lemma \ref{lem:d}, \eqref{eq:thmJcon}, \eqref{eq:thminv2}, \eqref{eq:thminv3} and \eqref{eq:Lie2inv22}, we obtain
\begin{eqnarray*}
dC_4^J&=&\frac{1}{2}d\huaS(\theta,J\circ\theta)\\
&=&\frac{1}{2}\Big(\huaS(d_{\nabla^0}\theta,J\circ\theta)-\huaS(\theta,d_{\nabla^1}(J\circ\theta))\Big)\\
&=&\frac{1}{2}\Big(\huaS(d_{\nabla^0}\theta,J\circ\theta)+\huaS(d_{\nabla^0}\theta,J\circ\theta) +2\huaS(R_\gamma,K\circ\theta)+2\huaS(\frkl_2\circ\theta,I_\gamma)\Big)\\
&=&\huaS(d_{\nabla^0}\theta,J\circ \theta)+\huaS(R_\gamma,K\circ\theta)+\huaS(\frkl_2\circ\theta,I_\gamma).
\end{eqnarray*}
The proof is finished. \qed\vspace{3mm}

For any $\huaG$-valued 1-form $\theta\in\Gamma(\Hom(TM,\huaG))$, we define a 4-form $C_4^{K}\in\Omega^4(M)$ by
\begin{equation}
  C_4^{K}=\frac{1}{3}\huaS(\theta,K\circ\theta).
\end{equation}
\begin{lem}\label{lem:thmK}
  With the above notations, we have
  \begin{equation}
    dC_4^K=\huaS(d_{\nabla^0}\theta,K\circ \theta)-\huaS(R_\gamma,\frkl_3\circ\theta)-\huaS(\frkl_2\circ\theta,J\circ \theta).
  \end{equation}
\end{lem}
\pf  By  Lemma \ref{lem:d}, \eqref{eq:thmKcon}, \eqref{eq:thminv2}, \eqref{eq:thminv3}, \eqref{eq:Lie2inv22} and \eqref{eq:Lie2inv33},  we have
\begin{eqnarray*}
   dC_4^K&=&\frac{1}{3}d\huaS(\theta,K\circ\theta)\\
   &=&\frac{1}{3}\Big(\huaS(d_{\nabla^0}\theta,K\circ\theta)-\huaS(\theta,d_{\nabla^1}(K\circ\theta))\Big)\\
    &=&\frac{1}{3}\Big(\huaS(d_{\nabla^0}\theta,K\circ\theta)+2\huaS(d_{\nabla^0}\theta,K\circ\theta) -3\huaS(R_\gamma,\frkl_3\circ\theta)-3\huaS(\frkl_2\circ\theta,J\circ \theta)\Big)\\
   & =&\huaS(d_{\nabla^0}\theta,K\circ \theta)-\huaS(R_\gamma,\frkl_3\circ\theta)-\huaS(\frkl_2\circ\theta,J\circ \theta).
\end{eqnarray*}
The proof is finished. \qed\vspace{3mm}

For any $\huaG$-valued 1-form $\theta\in\Gamma(\Hom(TM,\huaG))$, we define a 4-form $C_4^{\frkl_3}\in\Omega^4(M)$ by
\begin{equation}
  C_4^{\frkl_3}=\frac{1}{4}\huaS(\theta,\frkl_3\circ\theta).
\end{equation}
Equivalently,
$$
C_4^{\frkl_3}(X_1,X_2,X_3,X_4)=\huaS(\theta(X_1),\frkl_3(\theta(X_2),\theta(X_3),\theta(X_4))).
$$

\begin{lem}\label{lem:thml3}
 With the above notations, we have
 \begin{equation}
   dC_4^{\frkl_3}=\huaS(d_{\nabla^0}\theta,\frkl_3\circ\theta)+\huaS(\frkl_2\circ\theta,K\circ\theta).
 \end{equation}
\end{lem}
\pf By Lemma \ref{lem:d}, \eqref{eq:thmder3}, \eqref{eq:Lie2inv33} and \eqref{eq:thminv3}, we have
\begin{eqnarray*}
  dC_4^{\frkl_3}&=&\frac{1}{4}d\huaS(\theta,\frkl_3\circ\theta)\\
  &=&\frac{1}{4}\Big(\huaS(d_{\nabla^0}\theta,\frkl_3\circ\theta)-\huaS(\theta,d_{\nabla^1}\frkl_3\circ\theta)\Big)\\
  &=&\frac{1}{4}\Big(\huaS(d_{\nabla^0}\theta,\frkl_3\circ\theta)+3\huaS(d_{\nabla^0}\theta,\frkl_3\circ\theta)+4\huaS(\frkl_2\circ\theta,K\circ\theta)\Big)\\
  &=&\huaS(d_{\nabla^0}\theta,\frkl_3\circ\theta)+\huaS(\frkl_2\circ\theta,K\circ\theta).
\end{eqnarray*}
The proof is finished. \qed

\begin{lem}\label{lem:thml2l3}
  With the above notations, we have
  \begin{equation}
    \huaS(\frkl_2\circ\theta,\frkl_3\circ\theta)=0.
  \end{equation}
\end{lem}
\pf It follows from  \eqref{eq:Lie2inv22},  \eqref{eq:Lie2inv33} and the Jacobiator identity for $\frkl_3$. \qed

\section{Exact $\LWX$ 2-algebroids}
 A $\LWX$ 2-algebroid $(E_{-1},E_0,\partial,\rho,S,\diamond,\Omega)$ is called {\bf exact} if we have the following exact sequence of vector bundles:
$$
0\longrightarrow T^*M \stackrel{\rho^*}{\longrightarrow} E_{-1}\stackrel{\partial}{\longrightarrow}E_0\stackrel{\rho}{\longrightarrow}TM\longrightarrow0,
$$
where $\rho^*:T^*M\longrightarrow E_{-1}$ is defined by
$$
S(\rho^*(\alpha),e^0)= \langle\alpha,\rho(e^0)\rangle,\quad\forall \alpha\in\Omega^1(M), e^0\in\Gamma(E_0).
$$
\subsection{Skeletal exact $\LWX$ 2-algebroids}

In this subsection, we associate a skeletal exact $\LWX$ $2$-algebroid a cohomology class in $H^4(M)$, can call it the higher analogue of the  {\v{S}}evera class.

\begin{defi}
An exact $\LWX$ $2$-algebroid $(E_{-1},E_0,\partial,\rho,S,\diamond,\Omega)$ is called  {\bf skeletal} if $\partial=0.$
\end{defi}

For a skeletal exact $\LWX$ $2$-algebroid $(E_{-1},E_0,\partial,\rho,S,\diamond,\Omega)$, by the exactness, we have
$$
E_{-1}\cong T^*M,\quad E_0\cong TM.
$$
In the sequel, we will write $E_{-1}= T^*M,~E_0= TM$ directly. By definition, $\partial=0$. $S$ is given by
\begin{equation}\label{eq:Sske}
 S(X+\alpha,Y+\beta)=\langle\alpha,Y\rangle+\langle\beta,X\rangle,\quad \forall X,Y\in\frkX(M),\alpha,\beta\in\Omega^1(M).
\end{equation}
 Obviously,   we have
\begin{equation}\label{eq:brske1}
X\diamond Y=[X,Y].
\end{equation}
By Axiom (iv) in Definition \ref{defi:Courant-2 algebroid}, we have
\begin{eqnarray*}
  X\langle\alpha,Y\rangle&=&XS(\alpha,Y)=S(X\diamond\alpha,Y)+S(\alpha,X\diamond Y)\\
  &=&\langle X\diamond\alpha,Y\rangle+\langle\alpha,[X, Y]\rangle,
\end{eqnarray*}
which implies that
\begin{equation}\label{eq:brske2}
X\diamond\alpha=L_X\alpha.
\end{equation}
Then by Axiom (ii), we deduce that
\begin{equation}\label{eq:brske3}
  \alpha\diamond X=-\iota_Xd\alpha
\end{equation}
Finally, by Axiom (v), $\Omega:\wedge^3\frkX(M)\longrightarrow \Omega^1(M)$ gives rise to a $4$-form $H\in \Omega^4(M)$ by
\begin{equation}\label{eq:Hske}
  H(X,Y,Z,W)=S(\Omega(X,Y,Z),W).
\end{equation} By Axiom (i), $H$ is closed. Summarize the above discussion, we have
\begin{pro}
  Any skeletal exact $\LWX$ $2$-algebroid must be of the form $(T^*[1]M,TM,\partial=0,\rho={\id},S,\diamond,\Omega)$, where $S,\diamond$ are given by \eqref{eq:Sske}, \eqref{eq:brske1}, \eqref{eq:brske2}, \eqref{eq:brske3}, and $\Omega$ is equivalent to a closed $4$-form $H$.
\end{pro}

We will denote a skeletal exact $\LWX$ $2$-algebroid simply by $(T^*[1]M,TM,S,\diamond,H)$.

\begin{defi}
  An isomorphism from a $\LWX$ $2$-algebroid $(E_{-1},E_0,\partial,\rho,S,\diamond,\Omega)$ to a $\LWX$ $2$-algebroid $(E_{-1}',E_0',\partial',\rho',S',\diamond',\Omega')$ consists of two bundle isomorphisms $\Phi_0:E_0\longrightarrow E_0'$, $\Phi_{1}:E_{-1}\longrightarrow E_{-1}'$ and a bundle morphism $\Phi_2:\wedge^2E_0\longrightarrow E_{-1}'$ such that
  \begin{itemize}
    \item[{\rm (i)}] $\rho'\circ \Phi_0=\rho,$
     \item[{\rm (ii)}] $S(\Phi_0(e^0),\Phi_1(e^1))=S(e^0,e^1)$, for all $e^0\in\Gamma(E_0)$ and $e^1\in\Gamma(E_{-1})$,
     \item[{\rm (ii)}] $(\Phi_0,\Phi_1,\Phi_2)$ is an isomorphism from the Leibniz $2$-algebra $(\Gamma(E_{-1})\oplus\Gamma(E_0);\partial, \diamond,\Omega)$  to the Leibniz $2$-algebra $(\Gamma(E_{-1}')\oplus \Gamma(E_0'); \partial', \diamond', \Omega')$.
  \end{itemize}
\end{defi}

We have seen that different skeletal exact $\LWX$ 2-algebroids are differed by closed 4-forms. Now we show that if two closed 4-forms $H$ and $H'$ are in the same cohomology class, the corresponding skeletal exact $\LWX$ 2-algebroids are isomorphic. Thus, skeletal exact $\LWX$ 2-algebroids are classified by $H^4(M)$. This is a higher analogous result of that exact Courant algebroids are classified by $H^3(M)$.

\begin{thm}
 Let $(T^*[1]M,TM,S,\diamond,H)$ and $(T^*[1]M,TM,S,\diamond,H')$ be two skeletal exact $\LWX$ $2$-algebroids.
 They are isomorphic if and only if $H$ and $H'$ are in the same cohomology class.
\end{thm}
\pf If $H$ and $H'$ are in the same cohomology class, we assume that $H=H'+dh$. Then $(\Phi_0={\id},\Phi_1={\id},\Phi_2=h)$ gives the required isomorphism from  $(T^*[1]M,TM,S,\diamond,H)$ to $(T^*[1]M,TM,S,\diamond,H')$. In fact, what we need to show is that $(\Phi_0={\id},\Phi_1={\id},\Phi_2=h)$ is a Leibniz 2-algebra morphism from the Leibniz 2-algebra $(\Omega^1(M),\frkX(M),\partial=0,\diamond,H)$
to $(\Omega^1(M),\frkX(M),\partial=0,\diamond,H')$. The only nontrivial thing is to verify that
\begin{eqnarray*}
  X\diamond h(Y,Z)-Y\diamond h(X,Z)+h(X,Y)\diamond Z+H'(X,Y,Z)=h([X,Y],Z)+c.p.+H(X,Y,Z),
\end{eqnarray*}
which is equivalent to $dh+H'=H.$ \qed\vspace{3mm}

It is known that exact Courant algebroids are classified by the {\v{S}}evera class \cite{Severa:3-form}. Thus, for a skeletal exact $\LWX$ $2$-algebroid $(T^*[1]M,TM,S,\diamond,H)$, we call the cohomology class $[H]\in H^4(M)$ the {\bf higher analogue of the {\v{S}}evera class}.

\subsection{$\LWX$-extension of a quadratic Lie 2-algebroid}

In this subsection, we show that every exact $\LWX$ 2-algebroid gives rise to a quadratic Lie 2-algebroid and a quadratic Lie 2-algebroid admits a $\LWX$-extension if and only if its first Pontryagin class vanishes.

Let $(E_{-1},E_0,\partial,\rho,S,\diamond,\Omega)$ be an exact $\LWX$ 2-algebroid. Denote by $F_{-1}=E_{-1}/\rho^*(T^*[1]M)$ and we have the following short exact sequence:
$$
0\longrightarrow T^*[1]M \stackrel{\rho^*}{\longrightarrow} E_{-1}\stackrel{\pr}{\longrightarrow}F_{-1}\longrightarrow 0.
$$
On the graded vector bundle $F_{-1} \oplus E_0$, define $l_1:F_{-1} \longrightarrow E_0$ by
\begin{equation}
l_1(m)=\partial (\widetilde{m}),\quad\forall m\in\Gamma(F_{-1}), ~\widetilde{m}\in\Gamma(E_{-1})~\mbox{~such that~} \pr(\widetilde{m})=m;
\end{equation}
define $l_2$ by
\begin{eqnarray}
  l_2 (e^0_1, e^0_2)=e^0_1\diamond e^0_2,\quad
  l_2(e^0, m)=\pr(e^0\diamond \widetilde{m}),\quad
  l_2(m ,e^0) =\pr(\widetilde{m}\diamond e^0);
\end{eqnarray}
and define $l_3:\wedge^3\Gamma(E_0)\longrightarrow \Gamma(F_{-1})$ by
\begin{equation}
l_3(e^0_1,e^0_2,e^0_3)=\pr\Omega(e^0_1,e^0_2,e^0_3).
\end{equation}
\begin{pro}
 Let $(E_{-1},E_0,\partial,\rho,S,\diamond,\Omega)$ be an exact $\LWX$ $2$-algebroid. Then $(F_{-1}\oplus E_0; \rho, l_1, \\l_2, l_3)$ is a quadratic Lie $2$-algebroid, which is called the {\bf ample Lie $2$-algebroid}. Consequently, the exact $\LWX$ $2$-algebroid $(E_{-1},E_0,\partial,\rho,S,\diamond,\Omega)$ is an extension, called {\bf $\LWX$-extension}, of the quadratic Lie $2$-algebroid $(F_{-1}\oplus E_0;\rho,l_1,l_2,l_3)$ by  $T^*[1]M$, i.e.
 \begin{equation}\label{eq:ext-a}
\CD
  0 @>   >>  T^*[1]M @>\rho^*>> E_{-1}  @>\pr>> F_{-1} @>  >> 0 \\
  @. @V 0 VV @V  \partial  VV @V l_1 VV @.  \\
  0 @>   >> 0 @>  >> E_{0} @>   >> E_0@>  >>0.
\endCD
\end{equation}
\end{pro}
\pf By the exactness, we can deduce that $l_1$ is well-defined. By \cite[Lemma 3.7]{LiuSheng}, we can deduce that $l_2$ is also well-defined. Since the obstruction of the operation $\diamond$ being skew-symmetric is given in $\rho^*(T^*[1]M)$, it follows that $l_2$ is skew-symmetric. Since $\partial \circ \rho^*=0$, we obtain a Lie 2-algebra $(\Gamma(F_{-1})\oplus \Gamma(E_0);l_1,l_2,l_3)$. Therefore, $(F_{-1}\oplus E_0;\rho,l_1,l_2,l_3)$ is a transitive Lie $2$-algebroid. It is obvious that the pairing $S$ vanishes restricting on $\rho^*(T^*[1]M)\oplus \ker(\rho)$. Thus, it induces a nondegenerate graded symmetric pairing $\huaS$ on $F_{-1}\oplus \ker(\rho)$ by
\begin{equation}
  \huaS(u,m)=S(u,\widetilde{m}),\quad \forall m\in\Gamma(F_{-1}), ~\widetilde{m}\in\Gamma(E_{-1})~\mbox{~such that~} \pr(\widetilde{m})=m.
\end{equation}
It is straightforward to see that invariant conditions in Definition \ref{defi:quadraticLie2} are satisfied. Thus, the transitive Lie 2-algebroid $(F_{-1}\oplus E_0;\rho,l_1,l_2,l_3)$ is a quadratic Lie 2-algebroid.\qed\vspace{3mm}

Let $(A_{-1}\oplus A_0;\rho,l_1,l_2,l_3)$ be a quadratic Lie 2-algebroid, in which $l_1$ is injective.
Now we use the notations as in Subsection \ref{subsec:tranLie2}.
Denote by $\huaG=\ker(\rho)$. By the nondegeneracy of $\huaS$, we have $A_{-1}\cong \huaG^*$ and we will write $A_{-1}= \huaG^*[1]$ directly. Then $\huaS$ is simply given by  \eqref{defi:huaS}.
The transitive Lie 2-algebroid $(\huaG^*[1]\oplus A_0;\rho,l_1,l_2,l_3)$ fits the following exact sequence:
\begin{equation}\label{eq:ext-b}
\CD
  0 @>  >>  \huaG^*[1] @> >> \huaG^*[1]  @>  >>  0 @>  >> 0 \\
  @. @V\frkl_1 VV @V  l_1  VV @V  VV @.  \\
  0 @>>> \huaG=\ker(\rho) @>  >> A_{0} @> \rho  >> TM@> >>0.
\endCD
\end{equation}
We use $(\frkl_1,\frkl_2,\frkl_3)$ to denote the graded bundle of quadratic Lie 2-algebra structures on $\huaG^*[1]\oplus \huaG$. Since $l_1$ is injective, $\frkl_1$ is an isomorphism. Choose a splitting $(\sigma,\gamma)$ of the quadratic Lie $2$-algebroid $(\huaG^*[1]\oplus A_0;\rho,l_1,l_2,l_3)$, then $A_0\cong TM\oplus \huaG$ and $\rho$ is exactly the projection $\pr_{TM}$. Define $\nabla^0:\Gamma(TM)\times \Gamma(\huaG)\longrightarrow \Gamma(\huaG)$, $\nabla^1:\Gamma(TM)\times \Gamma(\huaG^*)\longrightarrow \Gamma(\huaG^*)$ and  $R\in\Omega^2(TM,\huaG)$ by \eqref{eq:definabla} and \eqref{eq:defiR} respectively. Define totally skewsymmetric bundle maps $I:\wedge^3TM\longrightarrow \huaG^*$, $J:\wedge^2TM\otimes\huaG\longrightarrow \huaG^*$ and $K:TM\otimes\wedge^2\huaG\longrightarrow \huaG^*$  by \eqref{defi:I}-\eqref{defi:K}  respectively. Then $I_\gamma,J,K$ induce totally skew-symmetric bundle maps $I^\flat_\gamma:\wedge^2TM\otimes\huaG\longrightarrow T^*M$, $J^\flat:TM\otimes\wedge^2\huaG\longrightarrow T^*M$ and $K^\flat:\wedge^3\huaG\longrightarrow T^*M$ by
\begin{eqnarray*}
  \langle I^\flat_\gamma(X,Y,u),Z\rangle&=&-\langle u,I_\gamma(X,Y,Z)\rangle,\\
  \langle J^\flat(X,u,v),Y\rangle&=&-\langle v,J(X,u,Y)\rangle,\\
  \langle K^\flat(u,v,w),X\rangle&=&-\langle w,K(X,u,v)\rangle,
\end{eqnarray*}
for all $X,Y,Z\in\frkX(M)$ and $u,v,w\in\Gamma(\huaG).$

\begin{thm}
 Let $(A_{-1}\oplus A_0;\rho,l_1,l_2,l_3)$ be a quadratic Lie $2$-algebroid, in which $l_1$ is injective. Then it admits a $\LWX$-extension if and only if its first Pontryagin class $[\huaS(R_\gamma,I_\gamma)]\in H^5(M)$ is trivial.
\end{thm}
\pf Assume that the quadratic Lie 2-algebroid $(A_{-1}\oplus A_0;\rho,l_1,l_2,l_3)$ admits a $\LWX$-extension. Using the same notations as in Section 4, let $(\huaG^*\oplus\huaG;\frkl_1,\frkl_2,\frkl_3)$ be the corresponding graded bundle of Lie 2-algebras. Let $(E_{-1},E_0,\partial,\rho,S,\diamond,\Omega)$ be an exact $\LWX$ $2$-algebroid whose ample Lie 2-algebroid is $(A_{-1}\oplus A_0;\rho,l_1,l_2,l_3)$. Then we have $E_0=A_0$.
We go on choosing a section $\lambda$ on $\huaG^*[1]$ that is orthogonal to $\sigma$, i.e. $S(\lambda(m),\sigma(X))=0$ for all $m\in\Gamma(\huaG^*[1])$ and $X\in\frkX(M)$. It turns out that $E_{-1}\cong \rho^*(T^*[1]M)\oplus \lambda(\huaG^*[1])$.  Since $\rho$ is surjective, $\rho^*$ is injective, we deduce that $E_{-1}\cong T^*[1]M\oplus \huaG^*[1]$ and the nondegenerate bilinear form is exactly given by
\begin{equation}
  S(X+u,\alpha+m)=\langle\alpha,X\rangle+\langle m,u\rangle,\quad \forall X\in\frkX(M),\alpha\in \Omega^1(M), u\in\Gamma(\huaG), m\in\Gamma(\huaG^*).
\end{equation}
Define $H\in\Omega^4(M)$ by
\begin{equation}
  H(X,Y,Z,W)=  S(\Omega(\sigma(X),\sigma(Y),\sigma(Z)),\sigma(W)).
\end{equation}

Transfer the $\LWX$ 2-algebroid structure to $(T^*[1]M\oplus\huaG^*[1])\oplus(TM\oplus \huaG)$, we have
\begin{eqnarray}\label{eq:formdiamond0}\left\{\begin{array}{rcl}
\partial(\alpha+m)&=&\frkl_1(m),\\
  (X+u)\diamond (Y+v)&=&[X,Y]+\nabla^0_Xv-\nabla^0_Yu-R_\gamma(X,Y)+\frkl_2(u,v),\\
  (X+u)\diamond (\alpha+m)&=&L_X\alpha+\nabla^1_Xm+Q(X,m)+\frkl_2(u,m)+P_1(u,m),\\
  (\alpha+m)\diamond (X+u)&=&-\iota_Xd\alpha-\nabla^1_Xm-Q(X,m)+\frkl_2(m,u)+P_2(m,u),\\
 \Omega(X+u,Y+v,Z+w)&=&H(X,Y,Z)-I_\gamma^\flat(X,Y,w)-I^\flat_\gamma(X,v,Z)-I_\gamma^\flat(u,Y,Z)\\
 &&-J^\flat(X,v,w)-J^\flat(u,Y,w)-J^\flat(u,v,Z)+K^\flat(u,v,w)\\
  &&-I_\gamma(X,Y,Z)-J(X,Y,w)-J(X,v,Z)-J(u,Y,Z)\\&&+K(X,v,w)+K(u,Y,w)+K(u,v,Z)+\frkl_3(u,v,w),
\end{array}\right.
\end{eqnarray}
where $Q:\frkX(M)\times \Gamma(\huaG^*)\longrightarrow \Omega^1(M)$ is defined by
\begin{equation}
  \langle Q(X,m),Y\rangle=\langle m,R_\gamma(X,Y)\rangle,
\end{equation}
and $P_1:\Gamma(\huaG)\times \Gamma(\huaG^*)\longrightarrow \Omega^1(M)$ and $P_2:\Gamma(\huaG^*)\times \Gamma(\huaG)\longrightarrow \Omega^1(M)$ are defined by
\begin{eqnarray}
  \langle P_1(u,m),X\rangle=\langle m,\nabla^0_Xu\rangle,\quad \langle P_2(m,u),X\rangle=\langle u,\nabla^1_Xm\rangle.
\end{eqnarray}

Then by the Jacobiator identity that $\Omega$ should satisfy, for all $W, X,Y,Z\in \frkX(M)$,   we have
 \begin{eqnarray*}
&& W\diamond \Omega(X,Y,Z)-X\diamond \Omega(W,Y,Z)+Y\diamond \Omega(W,X,Z)+\Omega(W,X,Y)\diamond Z\\
&&-\Omega(W\diamond X,Y,Z)+\Omega(W\diamond Y,X,Z)-\Omega(W\diamond Z,X,Y)\\
&&-\Omega(X\diamond Y,W,Z)+\Omega(X\diamond Z,W,Y)-\Omega(Y\diamond Z,W,X)\\
&=& W\diamond (H(X,Y,Z)-I(X,Y,Z))-X\diamond (H(W,Y,Z)-I(W,Y,Z))\\
&&+Y\diamond (H(W,X,Z)-I(W,X,Z))+(H(W,X,Y)-I(W,X,Y))\diamond Z\\
&&-\Omega([W, X]-R_\gamma(W,X),Y,Z)+\Omega([W, Y]-R_\gamma(W,Y),X,Z)-\Omega([W, Z]-R_\gamma(W,Z),X,Y)\\
&&-\Omega([X,Y]-R_\gamma(X,Y),W,Z)+\Omega([X, Z]-R_\gamma(X,Z),W,Y)-\Omega([Y, Z]-R_\gamma(Y,Z),W,X)\\
&=&L_WH(X,Y,Z)-L_XH(W,Y,Z)+L_YH(W,X,Z)-\iota_Zd(H(W,X,Y))\\
&&-H([W, X],Y,Z)+H([W, Y],X,Z)-H([W, Z],X,Y)\\
&&-H([X,Y],W,Z)+H([X, Z],W,Y)-H([Y, Z],W,X)\\
&&-\nabla^1_WI_\gamma(X,Y,Z)-Q(W,I_\gamma(X,Y,Z))+\nabla^1_XI_\gamma(W,Y,Z)+Q(X,I_\gamma(W,Y,Z))\\
&&-\nabla^1_YI_\gamma(W,X,Z)-Q(Y,I_\gamma(W,X,Z))+\nabla^1_ZI_\gamma(W,X,Y)+Q(Z,I_\gamma(W,X,Y))\\
&&+I_\gamma([W, X],Y,Z)-I_\gamma([W, Y],X,Z)+I_\gamma([W, Z],X,Y)\\
&&+I_\gamma([X,Y],W,Z)-I_\gamma([X, Z],W,Y)+I_\gamma([Y, Z],W,X)\\
&&-I^\flat_\gamma(R_\gamma(W,X),Y,Z)+I^\flat_\gamma(R_\gamma(W,Y),X,Z)-I^\flat_\gamma(R_\gamma(W,Z),X,Y)\\
&&-I^\flat_\gamma(R_\gamma(X,Y),W,Z)+I^\flat_\gamma(R_\gamma(X,Z),W,Y)-I^\flat_\gamma(R_\gamma(Y,Z),W,X)\\
&&-J(R_\gamma(W,X),Y,Z)+J(R_\gamma(W,Y),X,Z)-J(R_\gamma(W,Z),X,Y)\\
&&-J(R_\gamma(X,Y),W,Z)+J(R_\gamma(X,Z),W,Y)-J(R_\gamma(Y,Z),W,X)\\
&=&dH(W,X,Y,Z,\cdot)+\huaS(R_\gamma,I_\gamma)(W,X,Y,Z,\cdot)-(d_{\nabla^1}I_\gamma+J\circ R_\gamma)(W,X,Y,Z)\\
&=&0.
 \end{eqnarray*}
 In particular, we have $\huaS(R_\gamma,I_\gamma)+dH=0,$ which implies that the first Pontryagin class $[\huaS(R_\gamma,I_\gamma)]\in H^5(M)$ is trivial.

 Conversely, if the first Pontryagin class is trivial, then there exists a 4-form $H\in\Omega^4(M)$ such that   \begin{eqnarray}\label{eq:exact5form}
           \huaS(  R_\gamma,I_\gamma ) +dH=0.
  \end{eqnarray} On the graded bundle $(T^*[1]M\oplus\huaG^*[1])\oplus(TM\oplus \huaG)$,  define $\partial,\diamond,\Omega$ by \eqref{eq:formdiamond0}. We are going to show that $((T^*[1]M\oplus\huaG^*[1])\oplus(TM\oplus \huaG),\partial,\pr_{TM},S,\diamond,\Omega)$ is a $\LWX$ 2-algebroid.
  First we show that under conditions in Theorem \ref{thm:transtiveLie2} and by Proposition \ref{pro:quadraticLie2} and \eqref{eq:exact5form},
 $(\Omega^1(M)\oplus\Gamma(\huaG^*),\frkX(M)\oplus \Gamma(\huaG),\partial,\diamond,\Omega)$ is a Leibniz 2-algebra. This is the most intrinsic part in the proof. By \eqref{eq:compatl1} and \eqref{eq:Lie2inv11}, we can deduce that Axioms (a), (b), (c) in Definition \ref{defi:2leibniz} hold by straightforward computations. Note that the restriction of $\diamond$ on $\wedge^2\Gamma(TM\oplus\huaG)$ is the same as the one for the transitive Lie 2-algebroid given in \eqref{eq:strformula}, we deduce that Axiom (d) in Definition \ref{defi:2leibniz} holds by the fact $\partial|_{\Omega^1(M)}=0.$ By straightforward computation, we have
 \begin{eqnarray*}
  && (X+u)\diamond((Y+v)\diamond(\alpha+m))-((X+u)\diamond(Y+v))\diamond(\alpha+m)\\
  &&-(Y+v)\diamond((X+u)\diamond(\alpha+m))\\
  &=&l_2(X+u,l_2(Y+v,m))-l_2(l_2(X+u,Y+v),m)-l_2(Y+v,l_2(X+u,m))\\
  &&+L_X(L_Y\alpha+Q(Y,m)+P_1(v,m))+Q(X,\nabla^1_Ym+\frkl_2(v,m))+P_1(u,\nabla^1_Ym+\frkl_2(v,m))\\
  &&-L_{[X,Y]}\alpha-Q([X,Y],m)-P_1(\nabla^0_Xv-\nabla^0_Yu-R_\gamma(X,Y)+\frkl_2(u,v),m)\\
  &&-L_Y(L_X\alpha+Q(X,m)+P_1(u,m))-Q(Y,\nabla^1_Xm+\frkl_2(u,m))-P_1(v,\nabla^1_Xm+\frkl_2(u,m)),
 \end{eqnarray*}
 where $l_2$ is given by \eqref{eq:strformula}. By Theorem \ref{thm:transtiveLie2}, we have
 \begin{eqnarray*}
  &&l_2(X+u,l_2(Y+v,m))-l_2(l_2(X+u,Y+v),m)-l_2(Y+v,l_2(X+u,m))\\
  &=&l_3(X+u,Y+v,\frkl_1(m))\\
  &=&-J(X,Y,\frkl_1(m))+K(X,v,\frkl_1(m))+K(u,Y,\frkl_1(m))+\frkl_3(u,v,\frkl_1(m)).
 \end{eqnarray*}
 Thus, Axiom $(e_1)$ in Definition \ref{defi:2leibniz} holds if and only if
 \begin{eqnarray*}
  &&L_X(Q(Y,m)+P_1(v,m))+Q(X,\nabla^1_Ym+\frkl_2(v,m))+P_1(u,\nabla^1_Ym+\frkl_2(v,m))\\
  &&-Q([X,Y],m)-P_1(\nabla^0_Xv-\nabla^0_Yu-R_\gamma(X,Y)+\frkl_2(u,v),m)\\
  &&-L_Y(Q(X,m)+P_1(u,m))-Q(Y,\nabla^1_Xm+\frkl_2(u,m))-P_1(v,\nabla^1_Xm+\frkl_2(u,m))\\
  &=&-I^\flat_\gamma(X,Y,\frkl_1(m))-J^\flat(X,v,\frkl_1(m))-J^\flat(u,Y,\frkl_1(m))+K^\flat(u,v,\frkl_1(m)),
 \end{eqnarray*}
 which can be obtained by \eqref{eq:Lie2inv11}, \eqref{eq:Lie2inv22}, \eqref{eq:thmRcon}, \eqref{eq:thmder1},  \eqref{eq:thminv1} and \eqref{eq:thmobsmor1}. We omit the details.

 Similarly, we can show that Axiom $(e_2)$ and $(e_3)$ in Definition \ref{defi:2leibniz} hold.

 The last step to show that $((\Omega^1(M)\oplus\Gamma(\huaG^*))\oplus (\frkX(M)\oplus \Gamma(\huaG));\partial,\diamond,\Omega)$ is a Leibniz 2-algebra is to show the Jacobiator identity for $\Omega$. Roughly speaking, the Jacobiator identity for $\Omega$ is due to \eqref{eq:thmder3}, \eqref{eq:thmJcon}, \eqref{eq:thmIcon}, \eqref{eq:thmKcon} and \eqref{eq:exact5form}.   We leave the proof   to readers.
 This finishes the proof of $((\Omega^1(M)\oplus\Gamma(\huaG^*))\oplus (\frkX(M)\oplus \Gamma(\huaG));\partial,\diamond,\Omega)$ being a Leibniz 2-algebra.

 By \eqref{eq:thminv1}, we have
\begin{eqnarray*}
  (X+u)\diamond (\alpha+m)+(\alpha+m)\diamond (X+u)&=&L_X\alpha-\iota_Xd\alpha  +P_1(u,m) +P_2(m,u)\\
  &=&d(\langle X,\alpha\rangle+\langle u,m\rangle)\\
  &=&dS(X+u,\alpha+m),
\end{eqnarray*}
 which implies that Axiom (ii) in Definition \ref{defi:Courant-2 algebroid} holds.

 Finally, by Proposition \ref{pro:quadraticLie2}, we can deduce that Axioms (iii)-(iv) in Definition \ref{defi:Courant-2 algebroid} hold directly. The proof is finished.
 \qed

\section{The first Pontryagin class of a trivial principle 2-bundle with a $\Gamma$-connection}

\subsection{Strict Lie 2-groups and strict Lie 2-algebras}

A group is a monoid where every element has an inverse. A 2-group is
a monoidal category where every object has a weak inverse and every
morphism has an inverse. Denote the category of smooth manifolds and
smooth maps by $\Diff$, a (semistrict) Lie 2-group is  a 2-group
in $\DiffCat$, where  $\DiffCat$ is the 2-category consisting
of categories, functors, and natural transformations in $\Diff$.
For more details, see \cite{baez:2algebras,baez:2gp}. Here we only give the
definition of a strict Lie 2-group.
\begin{defi}
A strict Lie $2$-group is a Lie groupoid $C$
such that
\begin{itemize}
\item[\rm(a)]
The space of morphisms $C_1$ and the space of objects $C_0$ are Lie
groups.
\item[\rm(b)] The source and the target $s,t:C_1\longrightarrow C_0$, the identity assigning function $i:C_0\longrightarrow C_1$ and
 the composition $\circ:C_1\times_{C_0}C_1\longrightarrow C_1$ are all Lie group morphisms.
\end{itemize}
\end{defi}

In the following we will denote the composition $\circ$ in the Lie
groupoid structure by $\cdot_\ve:C_1\times_{C_0}C_1\longrightarrow
C_1$ and call it the vertical multiplication. Denote the Lie 2-group
multiplication by $\cdot_\hh:C\times C\longrightarrow C$ and call it
the horizontal multiplication.

It is well known that strict Lie 2-groups can be described by
crossed modules of Lie groups.

\begin{defi}
A crossed module of Lie groups is a quadruple $(H_1,H_0,\Psi,\Phi)$,
which we denote simply by $\HH$, where $H_1$ and $H_0$ are Lie
groups, $\Psi:H_1\longrightarrow H_0$ is a Lie group morphism, and
$\Phi:H_0\times H_1\longrightarrow H_1 $ is an action of $H_0$ on
$H_1$ preserving the Lie group structure of $H_1$ such that the Lie
group morphism $\Psi$ is $H_0$-equivariant:
\begin{equation}\label{cm g 1}
\Psi(\Phi_g(h))=g\Psi(h)g^{-1},\quad \forall ~g\in H_0,~h\in H_1,
\end{equation}
and $\Psi$ satisfies the so called Perffer identity:
\begin{equation}\label{cm g 2}
\Phi_{\Psi(h)}(h^\prime)=hh^\prime h^{-1},\quad\forall ~h,h^\prime\in
H_1.
\end{equation}
\end{defi}

\begin{thm}\label{thm:cm and slg}
There is a one-to-one correspondence between crossed modules of Lie
groups and strict Lie $2$-groups.
\end{thm}
 Roughly speaking, given a crossed
module $(H_1,H_0,\Psi,\Phi)$ of Lie groups, there is a strict Lie
2-group for which $C_0=H_0$ and $C_1= H_0\ltimes H_1$, the
semidirect product of $H_0$ and $H_1$. In this strict Lie 2-group,
the source and target maps $s,~t:C_1\longrightarrow C_0$ are given
by
$$
s(g,h)=g,\quad t(g,h)= t(h)\cdot g,
$$
the vertical multiplication $ \cdot_\ve$ is given by:
\begin{equation}\label{m v}
(g^\prime,h^\prime)\cdot_\ve(g,h) =(g, h^\prime\cdot h),\quad
\mbox{where} \quad g^\prime= t(h)\cdot g,
\end{equation}
the horizontal multiplication $\cdot_\mathrm{h}$ is given by
\begin{equation}\label{m h}
(g,h)\cdot_\mathrm{h} (g^\prime,h^\prime)=(g\cdot
g^\prime,h\cdot\Phi_g h^\prime).
\end{equation}

\begin{defi}
A crossed module of Lie algebras is a quadruple
$(\frkh_1,\frkh_0,\psi,\phi)$, which we denote by $\h$, where
$\frkh_1$ and $\frkh_0$ are Lie algebras,
$\psi:\frkh_1\longrightarrow\frkh_0$ is a Lie algebra morphism and
$\phi:\frkh_0\longrightarrow\Der(\frkh_1)$ is an action of Lie
algebra $\frkh_0$ on Lie algebra $\frkh_1$ as a derivation, such
that
$$
\psi(\phi_u(m))=[u,\psi(m)]_{\frkh_0},\quad \phi_{\psi(m)}(p)=[m,p]_{\frkh_1},\quad \forall u\in\frkh_0, m,p\in\frkh_1.
$$
\end{defi}

\begin{ex}
For any Lie algebra $\frkk$, the adjoint action $\ad$ is a Lie
algebra morphism from $\frkk$ to $\Der(\frkk)$. Then
$(\frkk,\Der(\frkk),\ad,\id)$ is a crossed module of Lie algebras.
\end{ex}

\begin{thm}\label{thm:dgla and cm}
There is a one-to-one correspondence between $2$-term DGLAs (strict Lie $2$-algebras) and
crossed modules of Lie algebras.
\end{thm}

In short, the formula for the correspondence can be given as follows: A 2-term DGLA $(V_{-1}\oplus V_0;\frkl_1,\frkl_2)$
gives rise to   a Lie algebra crossed module with $\frkh_1=V_{-1}$ and
$\frkh_0=V_0$, where the Lie brackets are given by:
\begin{eqnarray*}
~[m,p]_{\frkh_1}=\frkl_2(\frkl_1(m),p),\quad[u,v]_{\frkh_0}=\frkl_2(u,v), \quad \forall~m,p\in V_{-1},
~~u,v\in V_0,
\end{eqnarray*}
and $\psi=\frkl_1$, $\phi:\frkh_0\longrightarrow\Der(\frkh_1)$ is given by
$ \phi_u(m)=\frkl_2(u,m). $ The DGLA structure gives the Jacobi identity for
$[\cdot,\cdot]_{\frkh_1}$ and $[\cdot,\cdot]_{\frkh_0}$, and various other
conditions for crossed modules.

Conversely, a crossed module $(\frkh_1,\frkh_0,\psi,\phi)$ gives rise
to a 2-term DGLA with  $V_{-1}=\frkh_1$, $V_0=\frkh_0$, $\frkl_1=\psi$,
and $\frkl_2$ given by:
\begin{eqnarray*}
 ~ \frkl_2(u,v)\triangleq[u,v]_{\frkh_0},\quad\frkl_2(u,m)\triangleq\phi_u(m),\quad \forall
~u,v\in\frkh_0, m\in\frkh_1.
\end{eqnarray*}

\subsection{The transitive Lie 2-algebroid associated to a trivial principle 2-bundle with a $\Gamma$-connection }

First we review the notion of a principle 2-bundle for a strict Lie 2-group $\Gamma$ on the basis of \cite{waldorf:global,woc11}. Let $\Gamma$ be a strict Lie 2-group corresponding to the Lie group crossed module $(H_1,H_0,\Psi,\Phi)$. Let $(\frkh_1,\frkh_0,\psi,\phi)$ be the crossed module of Lie algebras corresponding to $(H_1,H_0,\Psi,\Phi)$, and $(\frkh_1\oplus \frkh_0;\frkl_1,\frkl_2)$ the associated strict Lie 2-algebra.
\begin{defi}
  A principle $\Gamma$-$2$-bundle over a differential manifold $M$ is a Lie groupoid $\huaP$, a surjective submersion functor $\pi:\huaP\longrightarrow M_{\mathrm{dis}}$, and a smooth right action $R$ of $\Gamma$ on $\huaP$ that preserves $\pi$, such that the smooth functor $$
  (\pr_1,R):\huaP\times \Gamma\longrightarrow \huaP\times_M\huaP
  $$
  is a weak equivalence, where $M_{\mathrm{dis}}$ is the Lie groupoid with objects $M$ and only identity morphisms.
\end{defi}

It is obvious that $M_{\mathrm{dis}}\times \Gamma$ is a principle $\Gamma$-$2$-bundle, which is called the trivial principle $\Gamma$-$2$-bundle.

A $\Gamma$-connection on $M$ is a pair $(A,B)$ consisting of an $\frkh_0$-valued 1-form $A\in\Omega^1(M,\frkh_0)$ and an $\frkh_1$-valued 2-form $B\in\Omega^2(M,\frkh_1)$. The curvature $\curv$ and the fake-curvature $\fcurv$ of a $\Gamma$-connection $(A,B)$ are defined by
\begin{eqnarray*}
  \curv(A,B)&\triangleq& dB+\frkl_2(A,B)\in\Omega^3(M,\frkh_1),\\
    \fcurv(A,B)&\triangleq& dA+\half\frkl_2(A,A)-\frkl_1 (B)\in\Omega^2(M,\frkh_0).
\end{eqnarray*}
Every $\Gamma$-connection can give rise to a connection on the trivial principle $\Gamma$-$2$-bundle $M_{\mathrm{dis}}\times \Gamma$. See \cite[Lemma 5.4.1]{waldorf:global} for more details.

Given a $\Gamma$-connection $(A,B)$, we can construct a transitive Lie 2-algebroid, which can be viewed as the infinitesimal of the trivial principle $\Gamma$-$2$-bundle with a connection, as follows: let $A_{-1}=M\times \frkh_1$, $A_0=TM\oplus (M\times \frkh_0)$, and define $l_1,l_2,l_3$ by
\begin{eqnarray}\label{eq:strformulastrict}\left\{\begin{array}{rcl}
  l_1(m)&=&\frkl_1(m),\\
  l_2(X+u,Y+v)&=&[X,Y]+L_Xv+\frkl_2(A(X),v)-L_Yu-\frkl_2(A(Y),u)\\&&+\fcurv(A,B)(X,Y)+\frkl_2(u,v),\\
  l_2(X+u,m)&=&L_Xm+\frkl_2(A(X),m)+\frkl_2(u,m),\\
  l_3(X+u,Y+v,Z+w)&=&-\curv(A,B)(X,Y,Z)+\frkl_2(B(X,Y),w)\\
  &&+\frkl_2(B(Y,Z),u)+\frkl_2(B(Z,X),v),
  \end{array}\right.
\end{eqnarray}
for all $X,Y,Z\in\frkX(M),u,v,w\in\Gamma(M\times \frkh_0)$ and $m\in\Gamma(M\times \frkh_1)$.

\begin{thm}\label{thm:quaLie2infi}
   Let $(A,B)$ be a $\Gamma$-connection on $M$. Then $((M\times \frkh_1)\oplus (TM\oplus (M\times \frkh_0)); \rho=\pr_{TM},l_1,l_2,l_3)$ is a transitive Lie $2$-algebroid, where $l_i$ are given by \eqref{eq:strformulastrict}.
\end{thm}
\pf Compare to  \eqref{eq:strformula}, we write
\begin{eqnarray*}
  \nabla^0_Xu&=&L_Xu+\frkl_2(A(X),u),\\
  \nabla^1_Xm&=&L_Xm+\frkl_2(A(X),m),\\
  R&=&-dA-\half\frkl_2(A,A),\\
  \gamma&=&B,\\
  J(X,Y,w)&=&-\frkl_2(B(X,Y),w),
\end{eqnarray*}
and $I=0,~K=0,~\frkl_3=0.$ Now we have
$$
I_\gamma=I+d_{\nabla^1}\gamma=d_{\nabla^1}B=dB+\frkl_2(A,B)=\curv(A,B).
$$Then $((M\times \frkh_1)\oplus (TM\oplus (M\times \frkh_0)); \rho=\pr_{TM},l_1,l_2,l_3)$ is a transitive Lie $2$-algebroid if and only if \eqref{eq:compatl1}-\eqref{eq:thmIcon} hold. Since $l_1$ is $\CWM$-linear, \eqref{eq:compatl1} holds naturally. Since $K=0$ and $\frkl_3=0$, we can deduce that \eqref{eq:thmder1}-\eqref{eq:thmder3} hold. It is also straightforward to deduce that \eqref{eq:thmobsmor1} and \eqref{eq:thmobsmor2} hold. By $d\frkl_2(A,A)=2\frkl_2(dA,A)$ and the Jacobi identity for $\frkl_2$, we have
\begin{eqnarray*}
  &&\nabla^0_XR(Y,Z)-R([X,Y],Z)+c.p.\\
  &=&L_XR(Y,Z)-R([X,Y],Z)+c.p.+\frkl_2(A(X),R(Y,Z))+c.p.\\
  &=&dR(X,Y,Z)+\frkl_2(A,R)(X,Y,Z)\\
  &=&d(-dA-\half\frkl_2(A,A))(X,Y,Z)+\frkl_2(A,-dA-\half\frkl_2(A,A))(X,Y,Z)\\
  &=&0,
\end{eqnarray*}
which implies that \eqref{eq:thmRcon} holds. It is straightforward to deduce that \eqref{eq:thmJcon} holds. Finally, by $d\frkl_2(A,\gamma)=\frkl_2(dA,\gamma)-\frkl_2(A,d\gamma)$ and $\frkl_2(\frkl_1(m),p)=\frkl_2(m,\frkl_1(p))$, we have
\begin{eqnarray*}
 d_{\nabla^1}(I+ d_{\nabla^1}\gamma)&=&d_{\nabla^1}(d_{\nabla^1}\gamma)\\
 &=&d(d_{\nabla^1}\gamma)+\frkl_2(A,d_{\nabla^1}\gamma)\\
 &=&d(d\gamma+\frkl_2(A,\gamma))+\frkl_2(A,d\gamma+\frkl_2(A,\gamma))\\
 &=&\frkl_2(dA,\gamma)+\frkl_2(A,\frkl_2(A,\gamma))\\
 &=&\frkl_2(dA+\half\frkl_2(A,A),\gamma),
\end{eqnarray*}
and
\begin{eqnarray*}
  J\circ R_\gamma=-\frkl_2(\gamma,R_\gamma)=-\frkl_2(\gamma,R+\frkl_1\circ\gamma)=\frkl_2(R,\gamma),
\end{eqnarray*}
which implies that \eqref{eq:thmIcon} holds. Thus, $((M\times \frkh_1)\oplus (TM\oplus (M\times \frkh_0)); \rho=\pr_{TM},l_1,l_2,l_3)$ is a transitive Lie $2$-algebroid. \qed

\begin{rmk}
  As proposed by the referee, it is natural to consider the case of nontrivial principle $2$-bundle. It is known that principle $2$-bundles are classified by the nonabelian differential cohomology. See \cite{NW,waldorf:parallel,woc11} for more details. A nonabelian cocycle is represented by an open cover $\{U_i\}$ of $M$, together with a collection of smooth maps $g_{ij}:U_{ij}\longrightarrow H_0$ and $a_{ijk}:U_{ijk}\longrightarrow H_1$ such that $g_{ik}=\Psi(a_{ijk})g_{ij}g_{jk}$ and $a_{ikl}\Phi_{g_{kl}}a_{ijk}=a_{ijl}a_{jkl}$. One can see that $g_{ij}$ does not satisfy the cocycle condition. Thus, the naive idea of gluing  the  local standard model given in Theorem \ref{thm:quaLie2infi} to obtain a global object does not work. We will study this interesting question in the future.
\end{rmk}

\subsection{The  first Pontryagin class of the quadratic Lie 2-algebroid associated to a $\Gamma$-connection}
In this subsection, we show that the  first Pontryagin class of the quadratic Lie 2-algebroid associated to an $\Gamma$-connection is trivial. First we give an interesting example of quadratic strict Lie 2-algebras.

\begin{ex}{\rm
  Let $(\g,[\cdot,\cdot]_\g,K)$ be a quadratic Lie algebra. That is, $K$ is a symmetric nondegenerate bilinear form on $\g$. Denote by $K^\sharp$ the induced map from $\g$ to $\g^*$, i.e.
  $$
  K^\sharp(u)(v)=K(u,v),\quad\forall u,v\in\g.
  $$
  Then $K^\sharp$ is an isomorphism. On the graded vector space $\g^*[1]\oplus\g,$ define $\frkl_1$ and $\frkl_2$ by
  \begin{eqnarray*}
    \frkl_1&=&(K^\sharp)^{-1},\\
    \frkl_2(u+\xi,v+\eta)&=&[u,v]_\g+\ad_u^*\eta-\ad_v^*\xi,\quad \forall u,v\in\g, \xi,\eta\in\g^*.
  \end{eqnarray*}
  Then $(\g^*[1]\oplus\g;\frkl_1,\frkl_2,\huaS)$ is a quadratic strict Lie $2$-algebra, where $\huaS$ is given by
  $$
  \huaS(u+\xi,v+\eta)=\langle\xi,v\rangle+\langle\eta,u\rangle.
  $$
  In fact, the only nontrivial part of proving $(\g^*[1]\oplus\g;\frkl_1,\frkl_2)$ to be a strict Lie $2$-algebra is to show the equality
  $$
  \frkl_2((K^\sharp)^{-1}(\xi),\eta)= \frkl_2(\xi,(K^\sharp)^{-1}(\eta)),\quad (K^\sharp)^{-1}\frkl_2(u,\eta)=\frkl_2(u, (K^\sharp)^{-1}(\eta)).
  $$
 Since $K^\sharp$ is an isomorphism, we can assume that $\xi=K^\sharp(u)$ and $\eta=K^\sharp(v)$. Since $K$ is invariant, we have
 \begin{eqnarray*}
    \langle\frkl_2((K^\sharp)^{-1}(\xi),\eta)- \frkl_2(\xi,(K^\sharp)^{-1}(\eta)),w\rangle&=&\langle\frkl_2(u,K^\sharp(v))-\frkl_2(K^\sharp(u),v),w\rangle\\
    &=&\langle\ad_u^*K^\sharp(v)+\ad_v^*K^\sharp(u),w\rangle\\
    &=&-K(v,[u,w]_\g) -K(u,[v,w]_\g)\\
    &=&0.
 \end{eqnarray*}
 Similarly, for any $\gamma=K^\sharp(w)$, we have
 \begin{eqnarray*}
  \langle (K^\sharp)^{-1}\frkl_2(u,\eta)-\frkl_2(u, (K^\sharp)^{-1}(\eta)),\gamma\rangle&=& \langle \frkl_2(u,\eta),w\rangle-\langle \frkl_2(u,v),K^\sharp(w)\rangle\\
  &=&-K(v,[u,w]_\g)-K(w,[u,v]_\g)\\
  &=&0.
 \end{eqnarray*}
 Thus, $(\g^*[1]\oplus\g;\frkl_1,\frkl_2)$ is a strict Lie $2$-algebra. It is obvious that $\huaS$ is invariant. Therefore, $(\g^*[1]\oplus\g;\frkl_1,\frkl_2,\huaS)$ is a quadratic strict Lie $2$-algebra. \qed
 }
\end{ex}

Now let $\Gamma$ be a strict Lie 2-group such that the corresponding strict Lie 2-algebra is a quadratic strict Lie 2-algebra $(\frkh^*[1]\oplus \frkh;\frkl_1,\frkl_2,\huaS)$. Let $(A,B)$ be a $\Gamma$-connection on $M$. Then the transitive  Lie 2-algebroid given in Theorem \ref{thm:quaLie2infi} is a quadratic Lie 2-algebroid naturally.
Consider its first Pontryagin class, we have
\begin{thm}
Let $\Gamma$ be a strict Lie $2$-group such that the corresponding strict Lie $2$-algebra is a quadratic strict Lie $2$-algebra, and $(A,B)$  a $\Gamma$-connection on $M$. Then the first Pontryagin class associated to the quadratic Lie $2$-algebroid given in Theorem \ref{thm:quaLie2infi}, which is represented by the $5$-cocycle $-\huaS(\fcurv(A,B),\curv(A,B))$, is trivial.
\end{thm}

\pf First we have
\begin{eqnarray*}
  \huaS(\fcurv(A,B),\curv(A,B))&=&\huaS(dA+\half\frkl_2(A,A)-\frkl_1(B),dB+\frkl_2(A,B))\\
  &=&\huaS(dA,dB)+\huaS(dA,\frkl_2(A,B))+\huaS(\half \frkl_2(A,A),dB)\\
  &&+\huaS(\half\frkl_2(A,A),\frkl_2(A,B))-\huaS(\frkl_1(B),dB)-\huaS(\frkl_1(B),\frkl_2(A,B)).
\end{eqnarray*}

By the Jacobi identity that $\frkl_2$ satisfies and \eqref{eq:Lie2inv22}, we have
$$
\huaS(\half\frkl_2(A,A),\frkl_2(A,B))=-\huaS(\frac{1}{6}\frkl_2(A,\frkl_2(A,A)),B)=0.
$$

By  \eqref{eq:Lie2inv22} and the property that $\frkl_2$ being skew-symmetric, we have
$$
\huaS(\frkl_1(B),\frkl_2(A,B))=-\huaS(A,\frkl_2(\frkl_1(B),B))=0.
$$

It is not hard to see that $d\frkl_1(B)=\frkl_1dB$. Thus, by \eqref{eq:Lie2inv11}, we have
$$
\huaS(\frkl_1(B),dB)=\half d\huaS(\frkl_1(B),B).
$$

Finally by \eqref{eq:Lie2inv22}, we have
\begin{eqnarray*}
d\huaS(\half\frkl_2(A,A),B)&=&\huaS(\half d\frkl_2(A,A),B)+\huaS(\half\frkl_2(A,A),dB)\\
&=&\huaS(\frkl_2(dA,A),B)+\huaS(\half\frkl_2(A,A),dB)\\
&=&\huaS(dA,\frkl_2(A,B))+\huaS(\half\frkl_2(A,A),dB).
\end{eqnarray*}

Therefore, we have
$$
\huaS(\fcurv(A,B),\curv(A,B))=d\Big(\huaS(A,dB)+\huaS(\half\frkl_2(A,A),B)-\half \huaS(\frkl_1(B),B)\Big),
$$
which implies that the first Pontryagin class is trivial.\qed\vspace{3mm}

At the end of this subsection, we analyze how does the primitive form of the first Pontryagin class behave under the gauge transformation. Recall from \cite[Section 5.4]{waldorf:global} that a {\bf gauge transformation} between $\Gamma$-connections $(A,B)$ and $(A',B')$ on $M$ is a pair $(g,\phi)$ consisting of a smooth map $g:M\longrightarrow H_0$ and a 1-form $\phi\in\Omega^1(M,\frkh_1)$ such that
\begin{eqnarray}
  A'&=&\Ad_gA-g^*\bar{\theta}-\frkl_1(\phi),\\
  B'&=&(\Phi_g)_*B-d\phi+\half\frkl_2(\frkl_1(\phi),\phi)-\frkl_2(\Ad_\g A,\phi)+\frkl_2(g^*\bar{\theta},\phi),
\end{eqnarray}
 where $\bar{\theta}$ is the Maurer-Cartan $1$-form on the Lie group $H_0$.

 \begin{pro}
   Let $(A,B)$ and $(A',B')$ be two gauge equivalent $\Gamma$-connections on $M$. Then we have
   \begin{equation}
     \huaS(\fcurv(A',B'),\curv(A',B'))=\huaS(\fcurv(A,B),\curv(A,B)).
   \end{equation}
 \end{pro}
 \pf By the Maurer-Cartan equation that $\bar{\theta}$ satisfies, we have
 \begin{eqnarray*}
   \fcurv(A',B')&=&dA'+\half\frkl_2(A',A')-\frkl_1(B')\\
   &=&d\big(\Ad_g A-g^*\bar{\theta}-\frkl_1(\phi)\big)+\half\frkl_2\big(\Ad_g A-g^*\bar{\theta}-\frkl_1(\phi),\Ad_g A-g^*\bar{\theta}-\frkl_1(\phi)\big)\\
   &&-\frkl_1\big((\Phi_g)_*B-d\phi+\half\frkl_2(\frkl_1(\phi),\phi)-\frkl_2(\Ad_g A,\phi)+\frkl_2(g^*\bar{\theta},\phi)\big)\\
   &=&\Ad_g dA+\frkl_2(g^*\bar{\theta},\Ad_gA)-dg^*\bar{\theta}-\frkl_1(d\phi)\\
   &&+\half\Ad_g\frkl_2(A,A)-\frkl_2(\Ad_gA,g^*\bar{\theta})-\frkl_2(\Ad_gA,\frkl_1(\phi))\\
   &&+\half\frkl_2(g^*\bar{\theta},g^*\bar{\theta})
   +\half\frkl_2(\frkl_1(\phi),\frkl_1(\phi))+\frkl_2(g^*\bar{\theta},\frkl_1(\phi))\\
   &&-\Ad_g\frkl_1(B)+\frkl_1(d\phi)-\half\frkl_2(\frkl_1(\phi),\frkl_1(\phi))+\frkl_2(\Ad_gA,\frkl_1(\phi))-\frkl_2(g^*\bar{\theta},\frkl_1(\phi))\\
   &=&\Ad_g(dA+\half\frkl_2(A,A)-\frkl_1(B))\\
   &=&\Ad_g\fcurv(A,B).
   \end{eqnarray*}
Similarly, by a tedious computation, we obtain
\begin{eqnarray*}
  \curv(A',B')=(\Phi_g)_*\curv(A,B)-\frkl_2(\Ad_g\fcurv(A,B),\phi).
\end{eqnarray*}
Therefore, by the invariance condition that $\huaS$ satisfies, we have
\begin{eqnarray*}
  \huaS(\fcurv(A',B'),\curv(A',B'))&=&\huaS(\Ad_g\fcurv(A,B),(\Phi_g)_*\curv(A,B)-\frkl_2(\Ad_g\fcurv(A,B),\phi))\\
 &=& \huaS(\fcurv(A,B),\curv(A,B))+\huaS(\Ad_g\frkl_2(\fcurv(A,B),\fcurv(A,B)),\phi)\\
 &=&\huaS(\fcurv(A,B),\curv(A,B)).
\end{eqnarray*}
The proof is finished.\qed

\end{document}